\documentstyle[12pt]{article}

\font\twelveof=msym10 at 12pt

\def\R{\mbox{\twelveof R}}
\def\C{\mbox{\twelveof C}}
\def\Z{\mbox{\twelveof Z}}
\def\N{\mbox{\twelveof N}}
\def\id{\mbox{\rm id}}
\def\cqfd{\qquad\qquad\vrule height 4pt depth 2pt width 5pt}
\def\case#1#2{{\textstyle{#1\over #2}}}
\def\tm{\tilde m}
\def\tiota{\tilde{\iota}}
\def\tDelta{\tilde{\Delta}}
\def\tepsilon{\tilde{\epsilon}}
\def\tS{\tilde S}
\def\tgamma{\tilde{\gamma}}
\def\tg{\tilde g}
\def\tone{\tilde 1}
\def\tH{\tilde H}
\def\th{\tilde h}
\def\hJ{\hat J}
\def\hZ{\hat Z}
\def\hd{\hat d}
\def\half{\frac{1}{2}}
\def\bX{\overline{X}}
\def\bx{\overline{x}}
\def\ss{\scriptstyle}

\renewcommand{\theequation}{\arabic{section}.\arabic{equation}}

\newcommand{\dotimes}{\stackrel{}{\dot{\otimes}}}
\newcommand{\ddotimes}{\stackrel{\nu}{\ddot{\otimes}}}
\newcommand{\nudot}{\stackrel{\nu}{\cdot}}
\newcommand{\bin}[2]{\left(\begin{array}{c}#1\\#2\end{array}\right)}
\newtheorem{definition}{Definition}[section]
\newtheorem{proposition}[definition]{Proposition}
\newtheorem{corollary}[definition]{Corollary}

\title{
\hspace{9cm}{\normalsize ULB/229/CQ/97/3}\\
\vspace{1.5cm}
Duals of coloured quantum universal enveloping algebras and coloured
universal $\cal T$-matrices}
\author{C. Quesne\thanks{Directeur de recherches FNRS; Electronic mail:
cquesne@ulb.ac.be} \\
{\small \sl Physique Nucl\'eaire Th\'eorique et Physique Math\'ematique,
Universit\'e Libre de Bruxelles,} \\
{\small \sl Campus de la Plaine CP229, Boulevard~du Triomphe, B-1050 Brussels,
Belgium}}
\date{ }
\begin{document}
\baselineskip=20pt plus 1pt minus 1pt
\maketitle

\begin{abstract}
We extend the notion of dually conjugate Hopf (super)algebras to the
coloured Hopf
(super)algebras~${\cal H}^c$ that we recently introduced. We show that if the
standard Hopf (super)algebras~${\cal H}_q$ that are the building blocks
of~${\cal
H}^c$ have Hopf duals~${\cal H}_q^*$, then the latter may be used to construct
coloured Hopf duals~${\cal H}^{c*}$, endowed with coloured algebra and antipode
maps, but with a standard coalgebraic structure. Next, we review the case where
the ${\cal H}_q$'s are quantum universal enveloping algebras of Lie
(super)algebras~$U_q(g)$, so that the corresponding ${\cal H}_q^*$'s are quantum
(super)groups~$G_q$. We extend the Fronsdal and Galindo universal ${\cal
T}$-matrix formalism to the coloured pairs $\left(U^c(g), G^c\right)$ by
defining
coloured universal ${\cal T}$-matrices. We then show that together with the
coloured universal $\cal R$-matrices previously introduced, the latter
provide an
algebraic formulation of the coloured $RTT$-relations, proposed by Basu-Mallick.
This establishes a link between the coloured extensions of Drinfeld-Jimbo and
Faddeev-Reshetikhin-Takhtajan pictures of quantum groups and quantum algebras.
Finally, we illustrate the construction of coloured pairs by giving some
explicit
results for the two-parameter deformations of $\bigl(U(gl(2)),
Gl(2)\bigr)$, and
$\bigl(U(gl(1/1)), Gl(1/1)\bigr)$.
\end{abstract}

\vspace{0.5cm}

\hspace*{0.3cm}
PACS: 02.10.Tq, 02.20.Sv, 03.65.Fd, 11.30.Na

\hspace*{0.3cm}
Running title: Duals of coloured quantum algebras

\hspace*{0.3cm}
To be published in J. Math. Phys.
\newpage
%
%
\section{INTRODUCTION} \label{sec:intro}
In a recent paper~\cite{cq96a} (henceforth referred to as~I and whose equations
will be subsequently quoted by their number preceded by~I), we did
introduce some
new algebraic structures ${\cal H}^c = \left({\cal H}, {\cal C}, {\cal
G}\right)$,
termed coloured Hopf algebras (see also Ref.~\cite{cq96b}). They are constructed
by starting from a standard Hopf algebra set ${\cal H} = \{\,{\cal H}_q
\mid q \in
{\cal Q}\,\}$, corresponding to some parameter set~$\cal Q$, and by generalizing
the definitions of coalgebra maps, and antipodes, by combining the standard ones
with the transformations of an algebra isomorphism group~$\cal G$, called colour
group. Such transformations are labelled by some colour parameters, taking
values
in a colour set~$\cal C$. Whenever the starting Hopf algebras~${\cal H}_q$ are
quasitriangular, the resulting coloured one~${\cal H}^c$ is characterized by the
existence of a coloured universal $\cal R$-matrix, denoted by~${\cal R}^c$, and
satisfying the coloured Yang-Baxter equation (YBE), i.e., the YBE with
nonadditive
spectral parameters~\cite{bazhanov,footnote1}.\par
%
%
In~I, we applied these new concepts to the Drinfeld-Jimbo formulation of quantum
algebras~\cite{drinfeld}, and constructed various examples of coloured quantum
universal enveloping algebras (QUEA's) of both semisimple and nonsemisimple Lie
algebras.\par
%
%
In a more recent paper~\cite{cq97} (henceforth referred to as~II), we
showed that
our definitions can be easily extended to deal with graded algebraic
structures~\cite{chaichian}, thereby leading to coloured Hopf superalgebras in
general, and to coloured QUEA's of Lie superalgebras in particular. In such
cases,
quasitriangularity implies that
${\cal R}^c$ is a solution of the coloured graded YBE.\par
%
%
In the literature, there exist other extensions of Hopf algebras related to the
coloured~YBE, and using either the Drinfeld-Jimbo~\cite{drinfeld}, or the
Faddeev-Reshetikhin-Takhtajan~\cite{faddeev} approach to quantum groups and
quantum algebras. In the former class, one should mention two previous
introductions of some algebraic structures called coloured Hopf algebras,
the first
one by Ohtsuki~\cite{ohtsuki} in the context of knot theory, and the second
one by
Bonatsos {\it et al}~\cite{bonatsos} in the study of some nonlinear deformation
of~$su(2)$. In the latter class, some extensions of the
Faddeev-Reshetikhin-Takhtajan $RLL$- and $RTT$-relations, where $R$ is replaced
by a coloured $R$-matrix, have been considered by Kundu and
Basu-Mallick~\cite{kundu,basu}, thereby leading to coloured $U(R)$ and
$A(R)$ Hopf
algebras, respectively.\par
%
%
As shown in~I, some direct connections exist between our coloured Hopf algebras,
and those previously introduced by Ohtsuki~\cite{ohtsuki}, and by Bonatsos
{\it et
al}~\cite{bonatsos}, as the former may be considered as generalizations of the
latter. On the contrary, possible relationships with coloured $U(R)$ and
$A(R)$ Hopf
algebras have not been investigated so far.\par
%
%
It is one of the aims of the present paper to fill in this gap for the
coloured $A(R)$
algebras, introduced by Basu-Mallick~\cite{basu}. For such a purpose, we shall
have to extend the notion of dually conjugate Hopf algebras~\cite{majid} to the
coloured context.\par
%
%
In the case of standard Hopf algebras, the interest of such a duality
concept was
recently highlighted  by Fronsdal and Galindo~\cite{fronsdal93}, who constructed
dual bases for the two-parameter deformations of~$Gl(2)$ and~$U(gl(2))$,
and used
them to build a universal $\cal T$-matrix, providing an algebraic formulation of
the $RTT$-relations, and yielding the quantum group generalization of the
familiar
exponential map from a Lie algebra to a Lie group. A similar construction was
carried out for some other quantum group and algebra pairs~\cite{fronsdal94,
bonechi,morozov,chakra96a}, and used to get representations of one member of the
pair from those of the other~\cite{finkelstein,jaga}. Moreover, the Fronsdal and
Galindo approach was also recently generalized to quantum supergroups and
superalgebras~\cite{chakra96b}.\par
%
%
We plan to show here that this universal $\cal T$-matrix formalism can be
extended to coloured QUEA's and their duals, and that the resulting coloured
universal $\cal T$-matrix, to be denoted by~${\cal T}^c$, provides the key
notion
for establishing a link between such duals and Basu-Mallick's coloured
$A(R)$~algebras.\par
%
%
This paper is organized as follows. In Sec.~\ref{sec:duals}, duals~${\cal
H}^{c*}$ of
coloured Hopf (super)algebras~${\cal H}^c$ are defined. In
Sec.~\ref{sec:QUEA}, the
case of coloured QUEA's of Lie (super)algebras is considered, and coloured
universal $\cal T$-matrices are introduced. In Secs.~\ref{sec:gl(2)}
and~\ref{sec:gl(1/1)}, the theory developed in the previous Section is
applied to
construct pairs $\left({\cal H}^c, {\cal H}^{c*}\right)$ for the two-parameter
deformations of $\bigl(U(gl(2)), Gl(2)\bigr)$, and $\bigl(U(gl(1/1)),
Gl(1/1)\bigr)$,
respectively. Section~\ref{sec:conclusion} contains some concluding remarks.\par
%
%
\section{DUALS OF COLOURED HOPF (SUPER)ALGE- BRAS}    \label{sec:duals}
\setcounter{equation}{0}
Let ${\cal H}^c = \left({\cal H}, {\cal C}, {\cal G}\right) = \left({\cal
H}_q, +, m_q,
\iota_q, \Delta^{\lambda,\mu}_{q,\nu}, \epsilon_{q,\nu}, S^{\mu}_{q,\nu}; k,
{\cal Q}, {\cal C}, {\cal G}\right)$ be a coloured Hopf algebra over some
field~$k$
($=\C$ or \R)~\cite{cq96a,cq96b}. Here $\cal Q$, $\cal C$, and ${\cal G} = \{\,
\sigma^{\nu}: {\cal H}_q \to {\cal H}_{q^{\nu}} \mid q, q^{\nu} \in {\cal
Q}, \nu \in
{\cal C}\,\}$ denote the parameter set, the colour set, and the colour group,
respectively. Hence, by definition, the $\sigma^{\nu}$'s satisfy
Eqs.~(I2.1)--(I2.4).
The maps $m_q: {\cal H}_q \otimes {\cal H}_q \to {\cal H}_q$, and $\iota_q:
k \to
{\cal H}_q$ are standard multiplication and unit maps, whereas
$\Delta^{\lambda,\mu}_{q,\nu} \equiv \left(\sigma^{\lambda} \otimes
\sigma^{\mu}\right) \circ \Delta_q \circ \sigma_{\nu}: {\cal H}_{q^{\nu}}
\to {\cal
H}_{q^{\lambda}} \otimes {\cal H}_{q^{\mu}}$, $\epsilon_{q,\nu} \equiv
\epsilon_q
\circ \sigma_{\nu}: {\cal H}_{q^{\nu}} \to k$, and $S^{\mu}_{q,\nu} \equiv
\sigma^{\mu} \circ S_q \circ \sigma_{\nu}: {\cal H}_{q^{\nu}} \to {\cal
H}_{q^{\mu}}$, defined in terms of standard comultiplication~$\Delta_q$,
counit~$\epsilon_q$, and antipode~$S_q$, with $\sigma_{\nu} \equiv
\left(\sigma^{\nu}\right)^{-1}$, are called coloured comultiplication, counit,
and antipode maps, respectively, and satisfy some generalized axioms, stated in
Proposition~II.2 of~I.\par
%
%
Let us now assume that the standard Hopf algebras~${\cal H}_q$, belonging
to~$\cal H$, have Hopf duals $\left({\cal H}_q^*, +, \tm_q, \tiota_q, \tDelta_q,
\tepsilon_q, \tS_q; k\right)$ (or in short ${\cal H}_q^*$), and let us set
${\cal H}^*
= \{\,{\cal H}_q^* \mid q \in {\cal Q}\,\}$. This means~\cite{majid} that
for any $q
\in {\cal Q}$, there exists a doubly nondegenerate bilinear form
$\langle\, , \,\rangle_q$, such that the bialgebra and antipode maps of the
dual pair
$\left({\cal H}_q, {\cal H}_q^*\right)$ are related by
\begin{eqnarray}
  \left\langle \tm_q(x \otimes y), X \right\rangle_q & = & \left\langle x
\otimes y,
         \Delta_q(X) \right\rangle_q, \qquad \left\langle \tiota_q(1), X
         \right\rangle_q = \left\langle \tilde 1_q, X \right\rangle_q =
\epsilon_q(X),
         \label{eq:dual1} \\
  \left\langle \tDelta_q(x), X \otimes Y \right\rangle_q & = & \left\langle
x, m_q(X
         \otimes Y)\right\rangle_q, \qquad \tepsilon_q(x) = \left\langle x,
\iota_q(1)
         \right\rangle_q = \left\langle x, 1_q \right\rangle_q,
\label{eq:dual2} \\
  \left\langle \tS_q(x), X \right\rangle_q & = & \left\langle x, S_q(X)
         \right\rangle_q, \label{eq:dual3}
\end{eqnarray}
for any $x$, $y \in {\cal H}_q^*$, and any $X$, $Y \in {\cal H}_q$.\par
%
%
Let us next introduce
%
\begin{definition}   \label{def-rho}
Let $\rho^{\nu}: {\cal H}_q^* \to {\cal H}_{q^{\nu}}^*$ be defined by
\begin{equation}
  \left\langle \rho^{\nu}(x), \sigma^{\nu}(X) \right\rangle_{q^{\nu}} =
\left\langle
  x, X \right\rangle_q,   \label{eq:rho}
\end{equation}
for any $x \in {\cal H}_q^*$, $X \in {\cal H}_q$, $q \in {\cal Q}$,
and~$\nu \in {\cal
C}$.
\end{definition}
%
From the properties of the $\sigma^{\nu}$'s and the nondegeneracy
of~$\langle \, ,
\, \rangle_q$, it is then straightforward to show
%
\begin{proposition}    \label{prop-rho}
The maps~$\rho^{\nu}$, defined in Eq.~(\ref{eq:rho}), are coalgebra
isomorphisms,
and the set $\{\, \rho^{\nu} \mid \nu \in {\cal C}\,\}$ is a group isomorphic
to~$\cal G$ (and denoted by the same symbol) with respect to the composition of
maps. In other words, the $\rho^{\nu}$'s are one-to-one, and satisfy the
properties
\begin{equation}
  \tDelta_{q^{\nu}} \circ \rho^{\nu} = \left(\rho^{\nu} \otimes
\rho^{\nu}\right)
  \circ \tDelta_q, \qquad \tepsilon_{q^{\nu}} \circ \rho^{\nu} = \tepsilon_q,
\end{equation}
and
\begin{eqnarray}
  \forall \nu, \nu' \in {\cal C}: \rho^{\nu' \circ \nu} & = & \rho^{\nu'} \circ
          \rho^{\nu}: {\cal H}_q^* \to {\cal H}_{q^{\nu' \circ \nu}}^*,
\nonumber \\
  \rho^{\nu^0} & = & \id: {\cal H}_q^* \to  {\cal H}_{q^{\nu^0}}^* = {\cal
H}_q^*,
          \nonumber \\
  \forall \nu \in {\cal C}: \rho^{\nu^i} & = &\rho_{\nu} \equiv
          \left(\rho^{\nu}\right)^{-1}: {\cal H}_{q^{\nu}}^* \to {\cal
H}_{q}^*,
\end{eqnarray}
where $\nu' \circ \nu$, $\nu^0$, and~$\nu^i$ have the same meaning as in
Eqs.~(I2.2)--(I2.4).
\end{proposition}
\par
%
%
Proceeding as for~$\cal H$ in~I, we may now combine the definitions of ${\cal
H}^*$, $\cal C$, and~$\cal G$ into
%
\begin{definition}   \label{def-colmaps}
The maps $\tm^{\nu}_{q,\lambda,\mu}: {\cal H}_{q^{\lambda}}^* \otimes {\cal
H}_{q^{\mu}}^* \to {\cal H}_{q^{\nu}}^*$, $\tiota^{\nu}_q: k \to {\cal
H}_{q^{\nu}}^*$,
and $\tS^{\nu}_{q,\mu}: {\cal H}_{q^{\mu}}^* \to {\cal H}_{q^{\nu}}^*$,
defined by
\begin{equation}
  \tm^{\nu}_{q,\lambda,\mu} \equiv \rho^{\nu} \circ \tm_q \circ
  \left(\rho_{\lambda} \otimes \rho_{\mu}\right), \qquad \tiota^{\nu}_q \equiv
  \rho^{\nu} \circ \tiota_q, \qquad \tS^{\nu}_{q,\mu} \equiv \rho^{\nu}
\circ \tS_q
  \circ \rho_{\mu},   \label{eq:col-maps}
\end{equation}
for any $q \in {\cal Q}$, and any $\lambda$, $\mu$, $\nu \in {\cal C}$, are
called
coloured multiplication, unit, and antipode, respectively.
\end{definition}
\par
%
%
By using the dual pairing (\ref{eq:dual1})--(\ref{eq:dual3}),
Definitions~\ref{def-rho} and~\ref{def-colmaps}, as well as Definition~II.1 and
Proposition~II.2 of~I, it is easy to establish
%
\begin{proposition}   \label{prop-colaxioms}
The coloured multiplication, unit, and antipode maps, defined in
Eq.~(\ref{eq:col-maps}), are dual to the coloured comultiplication, counit, and
antipode of~${\cal H}^c$, i.e.,
\begin{eqnarray}
  \forall x \in {\cal H}_{q^{\lambda}}^*, \forall y \in {\cal
H}_{q^{\mu}}^*, \forall X
         \in {\cal H}_{q^{\nu}}: \left\langle \tm^{\nu}_{q,\lambda,\mu}(x
\otimes y),
         X \right\rangle_{q^{\nu}} & = & \left\langle x \otimes y,
\Delta^{\lambda,\mu}
         _{q,\nu}(X) \right\rangle_{q^{\lambda},q^{\mu}}, \nonumber \\
  \forall X \in {\cal H}_{q^{\nu}}: \left\langle \tiota^{\nu}_q(1), X
         \right\rangle_{q^{\nu}} & = & \epsilon_{q,\nu}(X), \nonumber \\
  \forall x \in {\cal H}_{q^{\mu}}^*, \forall X \in {\cal H}_{q^{\nu}}:
\left\langle
         \tS^{\nu}_{q,\mu}(x), X \right\rangle_{q^{\nu}} & = & \left\langle
x, S^{\mu}
         _{q,\nu}(X) \right\rangle_{q^{\mu}}.    \label{eq:col-dual}
\end{eqnarray}
In addition, they transform under~$\cal G$ as
\begin{eqnarray}
  \tm^{\nu}_{q,\alpha,\beta} \circ \left(\rho^{\alpha}_{\lambda} \otimes
           \rho^{\beta}_{\mu}\right) & = & \tm^{\nu}_{q,\lambda,\mu} =
           \rho^{\nu}_{\gamma} \circ \tm^{\gamma}_{q,\lambda,\mu}, \nonumber \\
  \rho^{\nu}_{\alpha} \circ \tiota^{\alpha}_q & = & \tiota^{\nu}_q, \nonumber \\
  \tS^{\nu}_{q,\alpha} \circ \rho^{\alpha}_{\mu} & = & \tS^{\nu}_{q,\mu} =
           \rho^{\nu}_{\beta} \circ \tS^{\beta}_{q,\mu},
\end{eqnarray}
and satisfy generalized associativity, unit, and antipode axioms
\begin{eqnarray}
  \tm^{\nu}_{q,\lambda,\mu} \circ \left(\tm^{\lambda}_{q,\alpha,\beta} \otimes
           \rho^{\mu}_{\gamma}\right) & = & \tm^{\nu}_{q,\lambda',\mu'} \circ
           \left(\rho^{\lambda'}_{\alpha} \otimes
           \tm^{\mu'}_{q,\beta,\gamma}\right), \nonumber \\
  \tm^{\nu}_{q,\lambda,\mu} \circ \left(\tiota^{\lambda}_q \otimes
           \rho^{\mu}_{\alpha}\right) & = & \tm^{\nu}_{q,\lambda',\mu'} \circ
           \left(\rho^{\lambda'}_{\alpha} \otimes \tiota^{\mu'}_q\right) =
           \rho^{\nu}_{\alpha}, \nonumber \\
  \tm^{\nu}_{q,\lambda,\mu} \circ \left(\tS^{\lambda}_{q,\alpha} \otimes
           \rho^{\mu}_{\alpha}\right) \circ \tDelta_{q^{\alpha}} & = &
           \tm^{\nu}_{q,\lambda',\mu'} \circ \left(\rho^{\lambda'}_{\alpha}
\otimes
           \tS^{\mu'}_{q,\alpha}\right) \circ \tDelta_{q^{\alpha}} =
\tiota^{\nu}_q
           \circ \tepsilon_{q^{\alpha}},    \label{eq:gen-alg}
\end{eqnarray}
as well as generalized bialgebra axioms
\begin{eqnarray}
  \tDelta_{q^{\nu}} \circ \tm^{\nu}_{q,\lambda,\mu} & = &
           \left(\tm^{\nu}_{q,\lambda,\mu} \otimes
\tm^{\nu}_{q,\lambda,\mu}\right)
           \circ (\id \otimes \tau \otimes \id) \circ
\left(\tDelta_{q^{\lambda}}
           \otimes \tDelta_{q^{\mu}}\right), \nonumber \\
  \tDelta_{q^{\nu}} \circ \tiota^{\nu}_q & = & \tiota^{\nu}_q \otimes
           \tiota^{\nu}_q, \nonumber \\
  \tepsilon_{q^{\nu}} \circ \tm^{\nu}_{q,\lambda,\mu} & = & \tepsilon_{q^{\lambda}}
           \otimes \tepsilon_{q^{\mu}}, \nonumber \\
  \tepsilon_{q^{\nu}} \circ \tiota^{\nu}_q & = & \id,   \label{eq:gen-bialg}
\end{eqnarray}
where $\rho^{\lambda}_{\mu} \equiv \rho^{\lambda} \circ \rho_{\mu}$, $\tau$ is
the twist map, and no summation is implied over repeated indices.
\end{proposition}
%
{\it Remarks.} (1) As usual for Hopf algebras, the twist map is defined by
$\tau (x
\otimes y) = y \otimes x$ for any $x$,~$y$ belonging to appropriate spaces.
(2)~As
shown in Eqs.~(\ref{eq:gen-alg}) and~(I2.7), the coloured antipodes
$\tS^{\nu}_{q,\mu}$, and~$S^{\mu}_{q,\nu}$ satisfy distinct generalized
axioms. A
similar remark is valid for the generalized bialgebra structures given in
Eqs.~(\ref{eq:gen-bialg}) and~(I2.8), respectively.\par
%
%
{}From Proposition~\ref{prop-colaxioms}, one obtains
%
\begin{corollary}
For any $q \in \cal Q$, any $\nu \in \cal C$, and $q_{\nu} \equiv q^{\nu^i}$,
$\left({\cal H}_q^*, +, \tm^{\nu}_{q_{\nu},\nu,\nu}, \tiota^{\nu}_{q_{\nu}},
\tDelta_q, \tepsilon_q, \tS^{\nu}_{q_{\nu},\nu}; k\right)$ is a Hopf
algebra over $k$
with multiplication~$\tm^{\nu}_{q_{\nu},\nu,\nu}$,
unit~$\tiota^{\nu}_{q_{\nu}}$, and antipode~$\tS^{\nu}_{q_{\nu},\nu}$,
defined by
particularizing Eq.~(\ref{eq:col-maps}). Moreover, it is the Hopf dual of
the Hopf
algebra $\left({\cal H}_q, +, m_q, \iota_q, \Delta^{\nu,\nu}_{q_{\nu},\nu},
\epsilon_{q_{\nu},\nu}, S^{\nu}_{q_{\nu},\nu}; k\right)$, considered in
Corollary~II.3 of~I.
\end{corollary}
%
{\it Remark.\/} In particular, for $\nu = \nu^0$, we get back the original Hopf
structures of~${\cal H}_q^*$ and~${\cal H}_q$.\par
%
%
This result can be generalized as follows:
%
\begin{definition}   \label{def-coldual}
If the standard Hopf algebras~${\cal H}_q$ that are the building blocks of a
coloured Hopf algebra~${\cal H}^c$ have Hopf duals~${\cal H}_q^*$, then the set
${\cal H}^* = \{\,{\cal H}_q^* \mid q \in {\cal Q}\,\}$, endowed with coloured
multiplication, unit, and antipode maps $\tm^{\nu}_{q,\lambda,\mu}$,
$\tiota^{\nu}_q$, $\tS^{\nu}_{q,\mu}$, as defined in (\ref{eq:col-maps}),
is called
coloured Hopf dual of~${\cal H}^c$, and denoted by any one of the symbols
$\left({\cal H}_q^*, +, \tm^{\nu}_{q,\lambda,\mu}, \tiota^{\nu}_q, \tDelta_q,
\tepsilon_q, \tS^{\nu}_{q,\mu}; k, {\cal Q}, {\cal C}, {\cal G}\right)$,
$\left({\cal
H}^*, {\cal C}, {\cal G}\right)$, or ${\cal H}^{c*}$.
\end{definition}
\par
%
%
The coloured antipode~$\tS^{\nu}_{q,\mu}$ satisfies some additional properties,
which are again distinct from the corresponding ones of~$S^{\mu}_{q,\nu}$, given
in Proposition~II.5 of~I.
%
\begin{proposition}   \label{prop-antipode}
The coloured antipode~$\tS^{\nu}_{q,\mu}$ of a coloured Hopf dual~${\cal
H}^{c*}$
fulfils the relations
\begin{eqnarray}
  \tS^{\nu}_{q,\gamma} \circ \tm^{\gamma}_{q,\alpha,\beta} & = &
           \tm^{\nu}_{q,\mu,\lambda} \circ \tau \circ
\left(\tS^{\lambda}_{q,\alpha}
           \otimes \tS^{\mu}_{q,\beta}\right), \qquad \tS^{\nu}_{q,\mu} \circ
           \tiota^{\mu}_q = \tiota^{\nu}_q,  \nonumber \\
  \left(\tS^{\nu}_{q,\mu} \otimes \tS^{\nu}_{q,\mu}\right) \circ
\tDelta_{q^{\mu}}
           & = & \tau \circ \tDelta_{q^{\nu}} \circ \tS^{\nu}_{q,\mu}, \qquad
           \tepsilon_{q^{\nu}} \circ \tS^{\nu}_{q,\mu} = \tepsilon_{q^{\mu}}.
\end{eqnarray}
\end{proposition}
\par
%
%
The case where the ${\cal H}_q$'s are Hopf superalgebras $\left({\cal H}_q, +,
\gamma_q, m_q, \iota_q, \Delta_q, \epsilon_q, S_q;k\right)$~\cite{chaichian} can
be dealt with only some minor changes. Here $\gamma_q: {\cal H}_q \to {\cal
H}_q$
denotes their grading automorphism, i.e.,
\begin{equation}
  \gamma_q(X) = (-1)^{\deg X} X,
\end{equation}
for any homogeneous $X \in {\cal H}_q$, where $\deg X = 0$ or~1 according to
whether $X$ is even or odd. As all vector spaces are now graded, the tensor
product
and the twist map are such that
\begin{equation}
  (X \otimes Y) (Z \otimes T) = (-1)^{(\deg Y)(\deg Z)} XZ \otimes YT, \qquad
  \tau(X \otimes Y) = (-1)^{(\deg X)(\deg Y)} Y \otimes X,
\label{eq:graded-spaces}
\end{equation}
for any homogeneous $X$, $Y$, $Z$, $T \in {\cal H}_q$.\par
%
%
Then, according to Definition~3.1 of~II, a coloured Hopf superalgebra\linebreak
$\left({\cal H}_q, +, \gamma_q, m_q, \iota_q, \Delta^{\lambda,\mu}_{q,\nu},
\epsilon_{q,\nu}, S^{\mu}_{q,\nu}; k, {\cal Q}, {\cal C}, {\cal G}\right)$
is a coloured
Hopf algebra, whose building blocks are Hopf superalgebras, and whose
colour group
elements~$\sigma^{\nu}$ are superalgebra isomorphisms. In other words, the
$\sigma^{\nu}$'s satisfy Eqs.~(I2.1)--(I2.4), and in addition are
grade-preserving
maps, i.e.,
\begin{equation}
  \sigma^{\nu} \circ \gamma_q = \gamma_{q^{\nu}} \circ \sigma^{\nu},
\end{equation}
for any $q \in \cal Q$, and $\nu \in \cal C$. Such coloured Hopf superalgebras,
which will still be denoted by~${\cal H}^c$, fulfil similar properties as
ordinary
coloured Hopf algebras, provided Eq.~(\ref{eq:graded-spaces}) is taken into
account.
In particular, if the ${\cal H}_q$'s are quasitriangular, then the coloured
universal
$\cal R$-matrix ${\cal R}^c = \{\, {\cal R}^{\lambda,\mu}_q \mid q \in {\cal Q},
\lambda, \mu \in {\cal C} \,\}$ of $\left({\cal H}^c, {\cal R}^c\right)$
satisfies the
coloured graded~YBE,
\begin{equation}
  R^{\lambda,\mu}_{q,12} R^{\lambda,\nu}_{q,13} R^{\mu,\nu}_{q,23} =
  R^{\mu,\nu}_{q,23} R^{\lambda,\nu}_{q,13} R^{\lambda,\mu}_{q,12},
\end{equation}
which has the same form as the nongraded one, but wherein the first relation in
Eq.~(\ref{eq:graded-spaces}) is used to evaluate the products.\par
%
%
Turning now to dual structures, let us assume that the Hopf superalgebras~${\cal
H}_q$ have Hopf duals $\left({\cal H}_q^*, +, \tgamma_q, \tm_q, \tiota_q,
\tDelta_q, \tepsilon_q, \tS_q; k\right)$ with corresponding grading maps denoted
by~$\tgamma_q$, i.e.,
\begin{equation}
  \tgamma_q(x) = (-1)^{\deg x} x,
\end{equation}
for any homogeneous $x \in {\cal H}_q^*$, and $\deg x = 0$ or~1 for $x$ even or
odd, respectively. The doubly nondegenerate bilinear forms $\langle \, , \,
\rangle_q$ not only satisfy Eqs.~(\ref{eq:dual1})--(\ref{eq:dual3}), but
are also
consistent or, in other words, fulfil the relation
\begin{equation}
  \left\langle \tgamma_q(x), X \right\rangle_q = \left\langle x, \gamma_q(X)
  \right\rangle_q,
\end{equation}
for any homogeneous $x \in {\cal H}_q^*$, and $X \in {\cal H}_q$.\par
%
%
It is then straightforward to show
%
\begin{proposition}
In the case of Hopf superalgebras, the maps $\rho^{\nu}: {\cal H}_q^* \to {\cal
H}_{q^{\nu}}^*$, defined in Eq.~(\ref{eq:rho}), preserve the grading of the
${\cal
H}_q^*$'s, i.e.,
\begin{equation}
  \rho^{\nu} \circ \tgamma_q = \tgamma_{q^{\nu}} \circ \rho^{\nu},
\end{equation}
for any $q \in \cal Q$, $\nu \in \cal C$, in addition to satisfying
Proposition~\ref{prop-rho}. Moreover the coloured multiplication, unit, and
antipode maps $\tm^{\nu}_{q,\lambda,\mu}$, $\tiota^{\nu}_q$,
$\tS^{\nu}_{q,\mu}$,
defined in Eq.~(\ref{eq:col-maps}), satisfy Propositions~\ref{prop-colaxioms}
and~\ref{prop-antipode}, provided the counterpart of
Eq.~(\ref{eq:graded-spaces})
for the duals~${\cal H}_q^*$ is taken into account.
\end{proposition}
%
Finally, Definition~\ref{def-coldual} can be extended to coloured Hopf duals of
coloured Hopf superalgebras with ${\cal H}^{c*} = \left({\cal H}_q^*, +,
\tgamma_q,
\tm^{\nu}_{q,\lambda,\mu}, \tiota^{\nu}_q, \tDelta_q, \tepsilon_q,
\tS^{\nu}_{q,\mu}; k, {\cal Q}, {\cal C}, {\cal G}\right)$.\par
%
%
\section{DUALS OF COLOURED QUEA'S OF LIE (SUPER)ALGEBRAS}    \label{sec:QUEA}
\setcounter{equation}{0}
Whenever the Hopf (super)algebras~${\cal H}_q$ are QUEA's of Lie
(super)algebras~$U_q(g)$, their Hopf duals~${\cal H}_q^*$ (if they exist)
are Hopf
(super)algebras of quantized functions on the corresponding Lie (super)groups,
$Fun_q(G) = G_q$.\par
%
%
Let $\{X_A\}$ and $\{x^A\}$ be dual bases of~${\cal H}_q$ and~${\cal H}_q^*$,
respectively. Hence
\begin{equation}
  \left\langle x^A, X_B \right\rangle_q = \delta^A_B.    \label{eq:pairing}
\end{equation}
The universal $\cal T$-matrix of $G_q$~\cite{fronsdal93,bonechi} is defined
as the
element of $G_q \otimes U_q(g)$ given by
\begin{equation}
  {\cal T}_q = \sum_A x^A \otimes X_A.    \label{eq:T}
\end{equation}
\par
%
%
In terms of the bases $\{X_A\}$ and~$\{x^A\}$, the duality
relations~(\ref{eq:dual2}), and~(\ref{eq:col-dual}), between ${\cal H}^c =
U^c(g)$
and ${\cal H}^{c*} = G^c$ can be written as
\begin{eqnarray}
  m_q(X_A \otimes X_B) & = & \sum_C f_{AB}^C(q) X_C, \qquad 1_q = \iota_q(1) =
          \sum_A g^A(q) X_A, \nonumber \\
  \Delta^{\lambda,\mu}_{q,\nu}(X_A) & = & \sum_{BC} h_A^{BC}\left(q^{\lambda},
          q^{\mu}, q^{\nu}\right) X_B \otimes X_C, \qquad
\epsilon_{q,\nu}(X_A) =
          \tg_A\left(q^{\nu}\right), \nonumber \\
  S^{\mu}_{q,\nu}(X_A) & = & \sum_B s_A^B\left(q^{\mu}, q^{\nu}\right) X_B,
          \label{eq:dualbis1}
\end{eqnarray}
and
\begin{eqnarray}
  \tm^{\nu}_{q,\lambda,\mu}(x^B \otimes x^C) & = & \sum_A
          h^{BC}_A\left(q^{\lambda}, q^{\mu}, q^{\nu}\right) x^A, \qquad
          \tiota^{\nu}_q(1) = \sum_A \tg_A(q^{\nu}) x^A, \nonumber \\
  \tDelta_q(x^C) & = & \sum_{AB} f^C_{AB}(q)\, x^A \otimes x^B, \qquad
          \tepsilon_q(x^A) = g^A(q), \nonumber \\
  \tS^{\nu}_{q,\mu}(x^B) & = & \sum_A s^B_A\left(q^{\mu}, q^{\nu}\right) x^A,
          \label{eq:dualbis2}
\end{eqnarray}
where the structure constants for the multiplication (resp.~coloured
comultiplication) in~$U^c(g)$ become the structure constants for the
comultiplication (resp.~coloured multiplication) in~$G^c$, and similarly for the
unit and counit, or the antipodes. Note that here we denote by the same symbol
basis elements~$X_A$ (resp.~$x^A$) belonging to different members ${\cal
H}_{q^{\nu}} = U_{q^{\nu}}(g)$ (resp.~${\cal H}_{q^{\nu}}^* = G_{q^{\nu}}$)
of the Hopf
algebra set~$\cal H$ (resp.~${\cal H}^*$) in order not to overload notation by
adding an extra index referring to the corresponding algebra. This should
cause no
confusion since from the context, it is always clear to which Hopf algebra $X_A$
(resp.~$x^A$) belongs.\par
%
%
Let us now introduce
%
\begin{definition}   \label{def-colT}
The set ${\cal T}^c = \{\, {\cal T}^{\lambda}_q \mid q \in {\cal Q},
\lambda \in {\cal
C} \,\}$, where ${\cal T}^{\lambda}_q \equiv {\cal T}_{q^{\lambda}}$ is
defined in
Eq.~(\ref{eq:T}), is called the universal ${\cal T}$-matrix of~$G^c$.
\end{definition}
%
\begin{definition}   \label{def-dprod}
For any $q \in \cal Q$, and any $\lambda$, $\mu$, $\nu$, $\alpha$, $\beta
\in \cal
C$, let ${\cal U}^{\lambda,\mu}_q \in G_{q^{\lambda}} \otimes U_{q^{\mu}}(g)$,
${\cal A}^{\alpha,\mu}_q \dotimes {\cal B}^{\beta,\mu}_q \in G_{q^{\alpha}}
\otimes
G_{q^{\beta}} \otimes U_{q^{\mu}}(g)$, and ${\cal A}^{\lambda,\alpha}_q
\ddotimes
{\cal B}^{\mu,\beta}_q \in G_{q^{\nu}} \otimes U_{q^{\alpha}}(g) \otimes
U_{q^{\beta}}(g)$ be defined by
\begin{eqnarray}
  {\cal U}^{\lambda,\mu}_q & \equiv & \left(\id \otimes
S^{\mu}_{q,\lambda}\right)
           \left({\cal T}^{\lambda}_q\right), \nonumber \\
  {\cal A}^{\alpha,\mu}_q \dotimes {\cal B}^{\beta,\mu}_q & \equiv & \sum_{ABCD}
           a_A^C\, b_B^D\, x^A \otimes x^B \otimes m_{q^{\mu}}(X_C \otimes X_D),
           \nonumber \\
  {\cal A}^{\lambda,\alpha}_q \ddotimes {\cal B}^{\mu,\beta}_q & \equiv &
           \sum_{ABCD} a_A^C\, b_B^D\, \tm^{\nu}_{q,\lambda,\mu}(x^A
\otimes x^B)
           \otimes X_C \otimes X_D,      \label{eq:ddprod}
\end{eqnarray}
where
\begin{equation}
  {\cal A}^{\lambda,\mu}_q \equiv \sum_{AC} a_A^C\, x^A \otimes X_C, \qquad
  {\cal B}^{\lambda,\mu}_q \equiv \sum_{BD} b_B^D\, x^B \otimes X_D,
\end{equation}
with $a_A^C$, $b_B^D \in k$, denote any two elements of $G_{q^{\lambda}} \otimes
U_{q^{\mu}}(g)$.
\end{definition}
\par
%
%
{}From Definitions~\ref{def-colT}, \ref{def-dprod}, Eqs.~(\ref{eq:dualbis1}),
(\ref{eq:dualbis2}), and the generalized antipode axiom for~$\tS^{\nu}_{q,\mu}$,
given in Eq.~(\ref{eq:gen-alg}), it is straightforward to show
%
\begin{proposition}   \label{prop-colT}
The elements of the coloured universal $\cal T$-matrix of~$G^c$ satisfy the
relations
\begin{eqnarray}
  {\cal T}^{\lambda}_q \dotimes {\cal T}^{\lambda}_q & = &
          \left(\tDelta_{q^{\lambda}} \otimes \id\right) \left({\cal
T}^{\lambda}_q
          \right), \nonumber \\
  {\cal T}^{\lambda}_q \ddotimes {\cal T}^{\mu}_q & = &
          \left(\id \otimes \Delta^{\lambda,\mu}_{q,{\nu}}\right)
          \left({\cal T}^{\nu}_q\right), \nonumber \\
  {\cal U}^{\lambda,\mu}_q & = & \left(\tS^{\lambda}_{q,\mu} \otimes \id\right)
          \left({\cal T}^{\mu}_q\right), \nonumber \\
  \left(\tm^{\nu}_{q,\lambda,\mu} \otimes \id\right) \left({\cal T}^{\lambda}_q
          \dotimes {\cal U}^{\mu,\lambda}_q\right) & = &
          \left(\tm^{\nu}_{q,\mu',\lambda} \otimes \id\right)
          \left({\cal U}^{\mu',\lambda}_q \dotimes {\cal T}^{\lambda}_q\right) =
          \tiota^{\nu}_q(1) \otimes 1_{q^{\lambda}}, \nonumber \\
  \left(\id \otimes m_{q^{\lambda}}\right) \left({\cal T}^{\lambda}_q
          \ddotimes {\cal U}^{\mu,\lambda}_q\right) & = &
          \left(\id \otimes m_{q^{\lambda}}\right)
          \left({\cal U}^{\mu',\lambda}_q \ddotimes {\cal
T}^{\lambda}_q\right) =
          \tiota^{\nu}_q(1) \otimes 1_{q^{\lambda}}.    \label{eq:prop-T}
\end{eqnarray}
\end{proposition}
%
{\it Remark.} The first relation in Eq.~(\ref{eq:prop-T}) is a well-known
property of
the universal $\cal T$-matrix of $G_{q^{\lambda}}$~\cite{bonechi}, but the
remaining relations are new. Whenever the colour parameters reduce to that
characterizing the unit element of~$\cal G$, they go over into corresponding
relations for the universal $\cal T$-matrix, with ${\cal U}^{\lambda,\mu}_q$
becoming the element of $G_q \otimes U_q(g)$ denoted as~${\cal T}_q^{-1}$ in
Ref.~\cite{bonechi}.\par
%
%
With the purpose of generalizing the algebraic formulation of the
$RTT$-relations~\cite{bonechi}, let us introduce the coloured counterparts of
${\cal T}_{q1} \equiv \sum_A x^A \otimes (X_A \otimes 1_q)$, and ${\cal T}_{q2}
\equiv \sum_A x^A \otimes (1_q \otimes X_A)$:
%
\begin{definition}
Let ${\cal T}^{\lambda}_{q,1\mu} \in G_{q^{\lambda}} \otimes U_{q^{\lambda}}(g)
\otimes U_{q^{\mu}}(g)$, and ${\cal T}^{\mu}_{q,2\lambda} \in G_{q^{\mu}}
\otimes U_{q^{\lambda}}(g) \otimes U_{q^{\mu}}(g)$ be defined by
\begin{equation}
  {\cal T}^{\lambda}_{q,1\mu} \equiv \sum_A x^A \otimes (X_A \otimes
  1_{q^{\mu}}), \qquad   {\cal T}^{\mu}_{q,2\lambda} \equiv \sum_B x^B \otimes
  (1_{q^{\lambda}} \otimes X_B),   \label{eq:T1-T2}
\end{equation}
where ${\cal T}^{\lambda}_q = \sum_A x^A \otimes X_A$, and ${\cal T}^{\mu}_q
= \sum_B x^B \otimes X_B$ are the universal $\cal T$-matrices
of~$G_{q^{\lambda}}$ and~$G_{q^{\mu}}$, respectively.
\end{definition}
%
The searched for result is expressed in
%
\begin{proposition}   \label{prop-colRTT}
If the coloured Hopf (super)algebra~$U^c(g)$ has a coloured universal $\cal
R$-matrix, then the coloured universal $\cal T$-matrix of the coloured Hopf
dual~$G^c$ satisfies the relation
\begin{equation}
  \left(\tone_{q^{\nu}} \otimes {\cal R}^{\lambda,\mu}_q\right)
  \left({\cal T}^{\lambda}_{q,1\mu} \nudot {\cal T}^{\mu}_{q,2\lambda}\right) =
  \left({\cal T}^{\mu}_{q,2\lambda} \nudot {\cal T}^{\lambda}_{q,1\mu}\right)
  \left(\tone_{q^{\nu}} \otimes {\cal R}^{\lambda,\mu}_q\right),
\label{eq:colRTT}
\end{equation}
where $\mbox{}\nudot\mbox{}$ denotes a generalized multiplication such that
\begin{eqnarray}
  {\cal T}^{\lambda}_{q,1\mu} \nudot {\cal T}^{\mu}_{q,2\lambda} & \equiv &
         \sum_{AB} \tm^{\nu}_{q,\lambda,\mu}(x^A \otimes x^B) \otimes (X_A
\otimes
         X_B),    \label{eq:genmul1}\\
  {\cal T}^{\mu}_{q,2\lambda} \nudot {\cal T}^{\lambda}_{q,1\mu} & \equiv &
         \sum_{AB} \tm^{\nu}_{q,\mu,\lambda}(x^B \otimes x^A) \otimes (X_A
\otimes
         X_B).    \label{eq:genmul2}
\end{eqnarray}
\end{proposition}
%
{\it Proof.} Comparing Eq.~(\ref{eq:genmul1}) with the last relation in
Eq.~(\ref{eq:ddprod}), and taking Proposition~\ref{prop-colT} into account, we
successively get
\begin{equation}
  {\cal T}^{\lambda}_{q,1\mu} \nudot {\cal T}^{\mu}_{q,2\lambda} =
  {\cal T}^{\lambda}_q \ddotimes {\cal T}^{\mu}_q = \left(\id \otimes
  \Delta^{\lambda,\mu}_{q,\nu}\right) \left({\cal T}^{\nu}_q\right).
  \label{eq:interm}
\end{equation}
Moreover, from Eqs.~(\ref{eq:genmul1})--(\ref{eq:interm}), we also obtain
\begin{equation}
  {\cal T}^{\mu}_{q,2\lambda} \nudot {\cal T}^{\lambda}_{q,1\mu} = (\id \otimes
  \tau) \left({\cal T}^{\mu}_{q,1\lambda} \nudot {\cal
T}^{\lambda}_{q,2\mu}\right)
  = \left(\id \otimes \tau \circ \Delta^{\mu,\lambda}_{q,\nu}\right)
  \left({\cal T}^{\nu}_q\right).
\end{equation}
Equation~(\ref{eq:colRTT}) then directly follows from Eq.~(I2.15),
expressing that
$U^c(g)$ is almost cocommutative.\cqfd
\par
%
%
To explicitly show that Proposition~\ref{prop-colRTT} provides an algebraic
formulation of the coloured $RTT$-relations, introduced by
Basu-Mallick~\cite{basu}, one has to consider matrix representations of
Eq.~(\ref{eq:colRTT}). For such a purpose, let us denote by~$D^{(i)}_q$ matrix
representations of~$U_q(g)$ in some $n_i$-dimensional $k$-modules~$V^{(i)}_q$,
where index~$i$ distinguishes between inequivalent representations. We may now
represent the elements of the coloured universal $\cal R$- and $\cal T$-matrices
by
\begin{equation}
  R^{\lambda(i),\mu(j)}_q \equiv \left(D^{(i)}_{q^{\lambda}} \otimes
D^{(j)}_{q^{\mu}}
  \right) \left({\cal R}^{\lambda,\mu}_q\right), \qquad T^{\lambda(i)}_q \equiv
  \sum_A x^A D^{(i)}_{q^{\lambda}}(X_A),   \label{eq:colR-rep}
\end{equation}
which are $n_i n_j \times n_i n_j$ and $n_i \times n_i$ matrices, valued in~$k$
and~$G_{q^{\lambda}}$, respectively. Such definitions can be extended to
provide
representations of the operators in Eqs.~(\ref{eq:T1-T2}), (\ref{eq:genmul1}),
and~(\ref{eq:genmul2}), as $n_i n_j \times n_i n_j$ matrices valued in
$G_{q^{\lambda}}$, $G_{q^{\mu}}$, $G_{q^{\nu}}$, and~$G_{q^{\nu}}$,
respectively:
\begin{eqnarray}
  T^{\lambda(i)}_{q,1\mu(j)} & \equiv & \sum_A x^A
\left(D^{(i)}_{q^{\lambda}}(X_A)
         \otimes I^{(j)}\right), \nonumber \\
  T^{\mu(j)}_{q,2\lambda(i)} & \equiv & \sum_B x^B \left(I^{(i)} \otimes
         D^{(j)}_{q^{\mu}}(X_B)\right), \nonumber \\
  T^{\lambda(i)}_{q,1\mu(j)} \nudot T^{\mu(j)}_{q,2\lambda(i)} & \equiv &
         \sum_{AB} \tm^{\nu}_{q,\lambda,\mu}(x^A \otimes x^B)
         \left(D^{(i)}_{q^{\lambda}}(X_A) \otimes D^{(j)}_{q^{\mu}}(X_B)\right),
         \nonumber \\
  T^{\mu(j)}_{q,2\lambda(i)} \nudot T^{\lambda(i)}_{q,1\mu(j)} & \equiv &
         \sum_{AB} \tm^{\nu}_{q,\mu,\lambda}(x^B \otimes x^A)
         \left(D^{(i)}_{q^{\lambda}}(X_A) \otimes
D^{(j)}_{q^{\mu}}(X_B)\right).
         \label{eq:extT-rep}
\end{eqnarray}
Here $I^{(i)}$, and~$I^{(j)}$ denote $n_i \times n_i$, and $n_j \times n_j$ unit
matrices.\par
%
%
In the case of QUEA's of Lie superalgebras, since the $D^{(i)}_q$'s are graded
matrices, their tensor product has to be evaluated in conformity with
Eq.~(\ref{eq:graded-spaces}). According to the usual
convention~\cite{chaichian},
matrix multiplication of graded matrices is assumed to be the same as that of
nongraded ones, whereas the tensor product of two graded matrices is defined by
\begin{equation}
  (A \otimes B)_{ij,kl} \equiv (-1)^{\pi_k(\pi_j + \pi_l)} A_{i,k} B_{j,l}.
  \label{eq:tensorprod}
\end{equation}
Here $\pi_k$ denotes the $\Z_2$-grade of the $k$th row or column of matrix~$A$
or~$B$. In particular, if $A$ is homogeneous, its degree is given by $\deg
A = \pi_k
+ \pi_l$ for any row index~$k$, and any column index~$l$. Special cases of
Eq.~(\ref{eq:tensorprod}) are
\begin{equation}
  (A \otimes I)_{ij,kl} = A_{i,k} \delta_{j,l}, \qquad (I \otimes A)_{ij,kl} =
  (-1)^{\pi_i(\pi_j + \pi_l)} \delta_{i,k} A_{j,l}.   \label{eq:A1-A2}
\end{equation}
\par
%
%
{}From Proposition~\ref{prop-colRTT}, it is now straightforward to obtain
%
\begin{corollary}
In the representation $D^{(i)}_{q^{\lambda}} \otimes D^{(j)}_{q^{\mu}}$ of
$U_{q^{\lambda}}(g) \otimes U_{q^{\mu}}(g)$, Eq.~(\ref{eq:colRTT}) becomes
\begin{equation}
  R^{\lambda(i),\mu(j)}_q \left(T^{\lambda(i)}_{q,1\mu(j)} \nudot
  T^{\mu(j)}_{q,2\lambda(i)}\right) = \left(T^{\mu(j)}_{q,2\lambda(i)} \nudot
  T^{\lambda(i)}_{q,1\mu(j)}\right) R^{\lambda(i),\mu(j)}_q,
\label{eq:colRTT-rep}
\end{equation}
where the various matrices are defined in Eqs.~(\ref{eq:colR-rep}),
and~(\ref{eq:extT-rep}).
\end{corollary}
%
{\it Remark.\/} If we consider, in particular, the defining representations
of~$U_{q^{\lambda}}(g)$ and~$U_{q^{\mu}}(g)$, Equation~(\ref{eq:colRTT-rep}) can
be rewritten in a simplified notation as
\begin{equation}
  R^{\lambda,\mu}_q \left(T^{\lambda}_{q1} \nudot T^{\mu}_{q2}\right) =
  \left(T^{\mu}_{q2} \nudot T^{\lambda}_{q1}\right) R^{\lambda,\mu}_q.
  \label{eq:colRTT-defrep}
\end{equation}
This relation is formally identical with the coloured $RTT$-relations
\begin{equation}
  R^{\lambda,\mu}_q \, T^{\lambda}_{q1} \, T^{\mu}_{q2} = T^{\mu}_{q2} \,
  T^{\lambda}_{q1} \, R^{\lambda,\mu}_q,   \label{eq:colRTT-basu}
\end{equation}
introduced by Basu-Mallick~\cite {basu}. The dependence upon~$\nu$ in
Eq.~(\ref{eq:colRTT-defrep}), which is absent in
Eq.~(\ref{eq:colRTT-basu}), only
means that we may evaluate the latter in various dual Hopf
algebras~$G_{q^{\nu}}$,
with $\nu$ running over~$\cal C$.\par
%
%
In the next two Sections, we shall proceed to illustrate the new concepts and
results presented here on two simple, but nevertheless significant examples.\par
%
%
\section{THE COLOURED TWO-PARAMETER QUEA $U^c(gl(2))$ AND ITS DUAL
$Gl^c(2)$}
\label{sec:gl(2)}
\setcounter{equation}{0}
In~I, we constructed a coloured QUEA by starting from the two-parameter
deformation of $U(gl(2))$~\cite{schirrmacher,dobrev} in the Burd\'\i k and
Hellinger
formulation~\cite{burdik92a}. On the other hand, Fronsdal and
Galindo~\cite{fronsdal93} analyzed the duality relationship between the pair of
Hopf algebras $U_{p,q}(gl(2))$ and~$Gl_{p,q}(2)$, and constructed the universal
$\cal T$-matrix of the latter. It is therefore interesting to study how
this duality
picture can be extended to the coloured context. For such a purpose, it
will prove
convenient to use the approach to~$U_{p,q}(gl(2))$ considered in
Refs.~\cite{fronsdal93,jaga}, instead of that of Ref.~\cite{burdik92a}. The
resulting coloured QUEA will therefore differ from that constructed in~I.
One shoud
remember in this respect that there exist many ways of transforming a given Hopf
algebra into a coloured one, depending upon the set of generators and the colour
group that are used.\par
%
%
\subsection{The dual Hopf algebras $U_{p,q}(gl(2))$ and $Gl_{p,q}(2)$}
In the standard formulation~\cite{jaga,schirrmacher}, the quantum algebra
$U_{p,q}(gl(2))$ is defined as the algebra generated by the elements $\{J_0,
J_{\pm}, Z\}$, subject to the relations
\begin{equation}
  \left[J_0, J_{\pm}\right] = \pm J_{\pm}, \qquad \left[J_+, J_-\right] =
  \left[2J_0\right]_Q, \qquad \left[Z, J_0\right] = \left[Z, J_{\pm}\right] = 0,
  \label{eq:gl2-alg}
\end{equation}
where as usual
\begin{equation}
  [X]_t \equiv \frac{t^X - t^{-X}}{t - t^{-1}},   \label{eq:box}
\end{equation}
and the following combinations of $p$,~$q$ are considered:
\begin{equation}
  Q \equiv \sqrt{pq}, \qquad P \equiv \sqrt{p/q}.   \label{eq:QP}
\end{equation}
Whereas the algebraic relations~(\ref{eq:gl2-alg}) only depend upon the first of
these combined parameters, the coalgebraic structure depends upon both of
them.\par
%
%
To describe the duality relationship between $U_{p,q}(gl(2))$
and~$Gl_{p,q}(2)$, one
has to consider another set of $U_{p,q}(gl(2))$ generators $\{\hJ_0, \hJ_{\pm},
\hZ\}$, defined in terms of the first one by~\cite{jaga}
\begin{equation}
  \hJ_0 = J_0, \qquad \hJ_{\pm} = J_{\pm}\, Q^{\mp J_0 - \half}\, P^{Z - \half},
  \qquad \hZ = Z,
\end{equation}
and satisfying the commutation relations
\begin{equation}
  \left[\hJ_0, \hJ_{\pm}\right] = \pm \hJ_{\pm}, \qquad \left[\hJ_+,
\hJ_-\right] =
  P^{2\hZ - 1} \left[2\hJ_0\right]_Q, \qquad \left[\hZ, \hJ_0\right] =
\left[\hZ,
  \hJ_{\pm}\right] = 0.
  \label{eq:gl2-algbis}
\end{equation}
The remaining Hopf algebra maps are
\begin{eqnarray}
  \Delta_{p,q}\left(\hJ_0\right) & = & \hJ_0 \otimes 1 + 1 \otimes \hJ_0, \qquad
         \Delta_{p,q}\left(\hZ\right) = \hZ \otimes 1 + 1 \otimes \hZ,
\nonumber \\
  \Delta_{p,q}\left(\hJ_+\right) & = & \hJ_+ \otimes Q^{-2\hJ_0} P^{2\hZ} + 1
         \otimes \hJ_+, \qquad \Delta_{p,q}\left(\hJ_-\right) = \hJ_- \otimes 1
         + Q^{2\hJ_0} P^{2\hZ} \otimes \hJ_-, \nonumber \\
  \epsilon_{p,q}(X) & = & 0, \qquad X \in \{\hJ_0, \hJ_{\pm}, \hZ\},
\nonumber \\
  S_{p,q}\left(\hJ_0\right) & = & - \hJ_0, \qquad S_{p,q}\left(\hZ\right) =
- \hZ,
         \nonumber \\
  S_{p,q}\left(\hJ_+\right) & = & - \hJ_+ Q^{2\hJ_0} P^{-2\hZ}, \qquad
         S_{p,q}\left(\hJ_-\right) = - Q^{-2\hJ_0} P^{-2\hZ} \hJ_-,
\end{eqnarray}
and the universal $\cal R$-matrix can be written as~\cite{jaga,chakra94}
\begin{equation}
  {\cal R}_{p,q} = Q^{-2\hJ_0 \otimes \hJ_0}\, P^{2\left(\hZ \otimes \hJ_0
- \hJ_0
  \otimes \hZ\right)} \sum_{n=0}^{\infty} \frac{(1 - Q^2)^n P^n}{[n]_Q!}\,
  Q^{-n(n+1)/2}\, \hJ_+^n \otimes \hJ_-^n,   \label{eq:gl2-univR}
\end{equation}
where $[n]_Q! \equiv [n]_Q [n-1]_Q \ldots [1]_Q$ for $n \in \N^+$, and
$[0]_Q! \equiv
1$.\par
%
%
In the $2 \times 2$ defining representation~$D_{p,q}$ of $U_{p,q}(gl(2))$,
 given by
\begin{eqnarray}
  D_{p,q}(\hJ_0) & = & \case{1}{2} \left(\begin{array}{cc}
                1 & 0 \\ 0 & -1
                \end{array} \right), \qquad
  D_{p,q}(\hZ) = \case{1}{2} \left(\begin{array}{cc}
                1 & 0 \\ 0 & 1
                \end{array} \right), \nonumber \\
  D_{p,q}(\hJ_+) & = & \left(\begin{array}{cc}
                0 & 1 \\ 0 & 0
                \end{array} \right), \qquad
  D_{p,q}(\hJ_-) = \left(\begin{array}{cc}
                0 & 0 \\ 1 & 0
                \end{array} \right),      \label{eq:gl2-defrep}
\end{eqnarray}
the universal $\cal R$-matrix~(\ref{eq:gl2-univR}) is represented by the $4
\times
4$ matrix
\begin{equation}
  R_{p,q} \equiv \left(D_{p,q} \otimes D_{p,q}\right) \left({\cal
R}_{p,q}\right) =
  Q^{1/2} \left(\begin{array}{cccc}
                      Q^{-1} & 0        & 0              & 0 \\[0.1cm]
                      0        & P^{-1} & Q^{-1} - Q & 0 \\[0.1cm]
                      0        & 0        & P              & 0 \\[0.1cm]
                      0        & 0        & 0              & Q^{-1}
              \end{array} \right).   \label{eq:gl2-R}
\end{equation}
\par
%
%
{}For the corresponding quantum group~$Gl_{p,q}(2)$, the defining $T$-matrix is
specified\linebreak by~\cite{fronsdal93,jaga}
\begin{equation}
  T_{p,q} = \left(\begin{array}{cc}
                  a & b \\[0.1cm]
                  c & d
                 \end{array}\right),    \label{eq:Gl2-T}
\end{equation}
with the commutation relations
\begin{eqnarray}
  ab & = & qba, \qquad ac = pca, \qquad bd = pdb, \qquad cd = qdc, \nonumber \\
  bc & = & (p/q) cb, \qquad ad - da = \left(q - p^{-1}\right) bc,
\label{eq:Gl2-alg}
\end{eqnarray}
following from the $RTT$-relations corresponding to the
$R$-matrix~(\ref{eq:gl2-R}). The coalgebra maps are given by
\begin{equation}
  \tDelta_{p,q}(T_{p,q}) = T_{p,q} \dotimes T_{p,q} =
  \left(\begin{array}{cc}
           a \otimes a + b \otimes c & a \otimes b + b \otimes d \\[0.1cm]
           c \otimes a + d \otimes c & c \otimes b + d \otimes d
           \end{array}\right), \qquad
  \tepsilon_{p,q}(T_{p,q}) = I.    \label{eq:Gl2-coalg}
\end{equation}
\par
%
%
The quantum determinant, defined by
\begin{equation}
  {\cal D} \equiv ad - qbc = ad - pcb = da - p^{-1}bc = da - q^{-1}cb,
  \label{eq:Gl2-det}
\end{equation}
is a group-like element, i.e.,
\begin{equation}
  \tDelta_{p,q}({\cal D}) = {\cal D} \otimes {\cal D}, \qquad
\tepsilon_{p,q}({\cal D})
  = 1,
\end{equation}
but is not central (except for~$p=q$), since
\begin{equation}
  {\cal D}a = a{\cal D}, \qquad {\cal D}b = P^{-2} b{\cal D}, \qquad
  {\cal D}c = P^2 c{\cal D}, \qquad {\cal D}d = d{\cal D}. \label{eq:Gl2-detcom}
\end{equation}
With the assumption that $\cal D$ is invertible, $Gl_{p,q}(2)$ is endowed
with an
antipode map
\begin{equation}
  \tS_{p,q}(T_{p,q}) = (T_{p,q})^{-1} = {\cal D}^{-1}
  \left(\begin{array}{cc}
            d     & - p^{-1} b \\[0.1cm]
            -pc & a
           \end{array}\right) =
  \left(\begin{array}{cc}
            d     & - q^{-1} b \\[0.1cm]
            -qc & a
           \end{array}\right) {\cal D}^{-1},    \label{eq:Gl2-S}
\end{equation}
such that $\tS_{p,q}({\cal D}) = {\cal D}^{-1}$. This completes the Hopf
algebraic
structure of~$Gl_{p,q}(2)$.\par
%
%
{}Following Fronsdal and Galindo~\cite{fronsdal93}, one may consider a Gauss
decomposition of the $T$-matrix~(\ref{eq:Gl2-T}),
\begin{equation}
  T_{p,q} =
      \left(\begin{array}{cc}
                 1            & 0 \\[0.1cm]
                 \gamma & 1
               \end{array}\right)
      \left(\begin{array}{cc}
                 a & 0 \\[0.1cm]
                 0 & \hd
               \end{array}\right)
      \left(\begin{array}{cc}
                 1 & \beta \\[0.1cm]
                 0 & 1
               \end{array}\right),   \label{eq:Gl2-gauss}
\end{equation}
where $a$ and $\hd$ are assumed invertible, thence $\beta = a^{-1} b$,
$\gamma = c
a^{-1}$, $\hd = d - c a^{-1} b$. By assuming that the algebra can be
augmented with
the logarithms of~$a$ and~$\hd$, and by using the maps
\begin{equation}
  a = e^{\alpha}, \qquad \hd = e^{-\delta},   \label{eq:Gl2-alphadelta}
\end{equation}
and the reparametrization
\begin{equation}
  P = e^{\theta}, \qquad Q = e^{\varphi},
\end{equation}
the quantum determinant~(\ref{eq:Gl2-det}) is transformed into ${\cal D} =
\exp(\alpha - \delta)$, and the new variables $\{\alpha, \beta, \gamma,
\delta\}$
are seen to satisfy the commutation relations of a solvable Lie algebra
\begin{eqnarray}
  [\alpha, \beta] & = & (\varphi - \theta) \beta, \qquad [\alpha, \gamma] =
(\varphi
          + \theta) \gamma, \qquad [\alpha, \delta] = 0, \nonumber \\[0cm]
  [\delta, \beta] & = & (\varphi + \theta) \beta, \qquad [\delta, \gamma] =
(\varphi
          - \theta) \gamma, \qquad [\beta, \gamma] = 0, \label{eq:Gl2-Lie}
\end{eqnarray}
with a noncocommutatitve coproduct structure. As a consequence,
$Gl_{p,q}(2)$ can
be embedded into the enveloping algebra $U_{p,q}(\{\alpha, \beta, \gamma,
\delta\})$ of such a Lie algebra. As shown by Fronsdal and
Galindo~\cite{fronsdal93}, $U_{p,q}(\{\alpha, \beta, \gamma, \delta\})$ is
dual to
the QUEA~$U_{p,q}(gl(2))$. Following common practice, we shall mean the larger
enveloping algebra~$U_{p,q}(\{\alpha, \beta, \gamma, \delta\})$ whenever
speaking  of the duality between $U_{p,q}(gl(2))$ and~$Gl_{p,q}(2)$.\par
%
%
Dual bases of~$U_{p,q}(gl(2))$ and~$Gl_{p,q}(2)$ are given by~\cite{fronsdal93}
\begin{equation}
  X_A = \frac{Q^{a_1(a_1-1)/2} \hJ_-^{a_1}}{[a_1]_Q!} \frac{H^{a_2}}{a_2!}
  \frac{\tH^{a_3}}{a_3!} \frac{Q^{-a_4(a_4-1)/2} \hJ_+^{a_4}}{[a_4]_Q!},
  \label{eq:gl2-basis}
\end{equation}
and
\begin{equation}
  x^A = \gamma^{a_1} \alpha^{a_2} \delta^{a_3} \beta^{a_4},
\label{eq:Gl2-basis}
\end{equation}
respectively, where
\begin{equation}
  H \equiv \hJ_0 + \hZ, \qquad \tH \equiv \hJ_0 - \hZ,
\end{equation}
$A = (a_1, a_2, a_3, a_4)$, and $a_1$, $a_2$, $a_3$,~$a_4 \in \N$. From
Eqs.~(\ref{eq:T}), (\ref{eq:gl2-basis}), and~(\ref{eq:Gl2-basis}), the
universal $\cal
T$-matrix of~$Gl_{p,q}(2)$ can be written as
\begin{equation}
  {\cal T}_{p,q} = {\cal E}xp_{Q^{-2}}\left(\gamma \hJ_-\right)
\exp\left(\alpha H +
  \delta \tH\right) {\cal E}xp_{Q^2}\left(\beta \hJ_+\right),
\label{eq:Gl2-univT}
\end{equation}
in terms of the basic exponential function
\begin{equation}
  {\cal E}xp_{t^2}(z) = \sum_{n=0}^{\infty} \frac{t^{-n(n-1)/2}}{[n]_t!} z^t.
\end{equation}
\par
%
%
In the next Subsection, it will prove convenient to consider another set of dual
bases
\begin{eqnarray}
  \bX_A & = & \sum_B c_A^B X_B = \frac{Q^{a_1(a_1-1)/2} \hJ_-^{a_1}}{[a_1]_Q!}
           \frac{\hJ_0^{a_2}}{a_2!} \frac{\hZ^{a_3}}{a_3!}
\frac{Q^{-a_4(a_4-1)/2}
           \hJ_+^{a_4}}{[a_4]_Q!}, \label{eq:gl2-basisbis} \\
  \bx^A & = & \sum_B \left(c^{-1}\right)^A_B x^B = \gamma^{a_1} h^{a_2}
\th^{a_3}
           \beta^{a_4}.  \label{eq:Gl2-basisbis}
\end{eqnarray}
Here, $A$ has the same meaning as in
Eqs.~(\ref{eq:gl2-basis}),~(\ref{eq:Gl2-basis}),
\begin{equation}
  h \equiv \alpha + \delta, \qquad \th \equiv \alpha - \delta,
\label{eq:Gl2-h}
\end{equation}
and
\begin{equation}
  c_A^B = 2^{-a_2-a_3} \left(c^{-1}\right)_A^B = \delta_{a_1}^{b_1}\,
  \delta_{a_2+a_3}^{b_2+b_3}\, \delta_{a_4}^{b_4}\, 2^{-a_2-a_3}
  \sum_{t=\max(0,b_2-a_2)}^{\min(b_2,a_3)} (-1)^{a_3-t} \bin{b_2}{t}
  \bin{b_3}{a_3-t},   \label{eq:basis-basisbis}
\end{equation}
where $\bin{s}{t}$ denotes a standard binomial coefficient. In terms of such
dual bases, the universal $\cal T$-matrix~(\ref{eq:Gl2-univT}) can be written as
\begin{equation}
  {\cal T}_{p,q} = {\cal E}xp_{Q^{-2}}\left(\gamma \hJ_-\right) \exp\left(h
\hJ_0 +
  \th \hZ\right) {\cal E}xp_{Q^2}\left(\beta \hJ_+\right).
\label{eq:Gl2-univTbis}
\end{equation}
\par
%
%
\subsection{The coloured Hopf algebra $U^c(gl(2))$ and its dual $Gl^c(2)$}
The defining relations~(\ref{eq:gl2-algbis}) of the
$U_{p,q}(gl(2))$~algebra are left
invariant under the transformations
\begin{equation}
  \sigma^{\nu}\left(\hJ_0\right) = \hJ_0, \qquad
  \sigma^{\nu}\left(\hJ_{\pm}\right) = P^{(\nu-1)/2} \hJ_{\pm}, \qquad
  \sigma^{\nu}\left(\hZ\right) = \nu \hZ,   \label{eq:gl2-sigma}
\end{equation}
where $\nu \in {\cal C} = \C \setminus \{0\}$, provided $P$ is replaced by its
$\nu$th power, while $Q$ is left unchanged, or equivalently~\cite{footnote2}
\begin{equation}
  (p,q) \to \left(p^{(\nu)}, q^{(\nu)}\right) = \left(p^{(1+\nu)/2}
q^{(1-\nu)/2},
  p^{(1-\nu)/2} q^{(1+\nu)/2}\right).  \label{eq:gl2-pqchange}
\end{equation}
Hence, the $\sigma^{\nu}$'s are isomorphic mappings between two
$U_{p,q}(gl(2))$~algebras with different parameters, as given in
Eq.~(\ref{eq:gl2-pqchange}), and they define a colour group ${\cal G} =
Gl(1,\C)$,
since $\nu' \circ \nu = \nu' \nu$, $\nu^0 = 1$, and~$\nu^i = \nu^{-1}$.\par
%
%
By choosing the parameter set ${\cal Q} = (\C \setminus \{0\}) \times (\C
\setminus \{0\})$, we therefore obtain a coloured quasitriangular Hopf
algebra~$U^c(gl(2))$, whose coloured comultiplication, counit, antipode,
and $\cal
R$-matrix are given by
\begin{eqnarray}
  \Delta^{\lambda,\mu}_{p,q,\nu}\left(\hJ_0\right) & = & \hJ_0 \otimes 1 + 1
           \otimes \hJ_0, \qquad
\Delta^{\lambda,\mu}_{p,q,\nu}\left(\hZ\right) =
           \frac{\lambda}{\nu}\, \hZ \otimes 1 + \frac{\mu}{\nu}\, 1
\otimes \hZ,
           \nonumber \\
  \Delta^{\lambda,\mu}_{p,q,\nu}\left(\hJ_+\right) & = & P^{(\lambda-\nu)/2}
           \hJ_+ \otimes Q^{-2\hJ_0} P^{2\mu\hZ} + P^{(\mu-\nu)/2} 1 \otimes
           \hJ_+, \nonumber \\
  \Delta^{\lambda,\mu}_{p,q,\nu}\left(\hJ_-\right) & = & P^{(\lambda-\nu)/2}
           \hJ_- \otimes 1 + P^{(\mu-\nu)/2} Q^{2\hJ_0} P^{2\lambda\hZ} \otimes
           \hJ_-, \nonumber \\
  \epsilon_{p,q,\nu}(X) & = & 0, \qquad X \in \{\hJ_0, \hJ_{\pm}, \hZ\},
\nonumber \\
  S^{\mu}_{p,q,\nu}\left(\hJ_0\right) & = & - \hJ_0, \qquad
           S^{\mu}_{p,q,\nu}\left(\hJ_+\right) = - P^{(\mu-\nu)/2} \hJ_+
Q^{2\hJ_0}
           P^{-2\mu\hZ}, \nonumber \\
  S^{\mu}_{p,q,\nu}\left(\hJ_-\right) & = & - P^{(\mu-\nu)/2} Q^{-2\hJ_0}
           P^{-2\mu\hZ} \hJ_-, \qquad S^{\mu}_{p,q,\nu}\left(\hZ\right) = -
           \frac{\mu}{\nu} \hZ, \nonumber \\
  {\cal R}^{\lambda,\mu}_{p,q} & = & Q^{-2 \hJ_0 \otimes \hJ_0} \,
           P^{2\left(\lambda \hZ \otimes \hJ_0 - \mu \hJ_0 \otimes \hZ\right)}
           \nonumber \\
  & & \mbox{} \times  \sum_{n=0}^{\infty} \frac{(1 - Q^2)^n}{[n]_Q!}
Q^{-n(n+1)/2}
           P^{(\lambda+\mu)n/2} \hJ_+^n \otimes \hJ_-^n,
\end{eqnarray}
respectively.\par
%
%
In the defining representation~(\ref{eq:gl2-defrep}), the coloured
universal $\cal
R$-matrix is represented by the $4 \times 4$~matrix
\begin{equation}
  R^{\lambda,\mu}_{p,q} = Q^{1/2}
     \left(\begin{array}{cccc}
          Q^{-1} P^{(\lambda-\mu)/2} & 0 & 0 & 0 \\[0.1cm]
          0 & P^{-(\lambda+\mu)/2} & Q^{-1} - Q & 0 \\[0.1cm]
          0 & 0 & P^{(\lambda+\mu)/2} & 0 \\[0.1cm]
          0 & 0 & 0 & Q^{-1} P^{-(\lambda-\mu)/2}
     \end{array} \right),   \label{eq:gl2-colR}
\end{equation}
which is a matrix solution of the coloured YBE~(I1.5), and gives back
matrix~(\ref{eq:gl2-R}) whenever $\lambda$,~$\mu\to 1$.\par
%
%
Under transformation~(\ref{eq:gl2-sigma}), the basis elements~$\bX_A$
of~$U_{p,q}(gl(2))$, defined in Eq.~(\ref{eq:gl2-basisbis}), become
\begin{equation}
  \sigma^{\nu}\left(\bX_A\right) = P^{(\nu-1)(a_1+a_4)/2} \nu^{a_3} \bX_A.
\end{equation}
Hence, from Eqs.~(\ref{eq:rho}) and~(\ref{eq:pairing}), it follows that
\begin{equation}
  \rho^{\nu}\left(\bx^A\right) = P^{(1-\nu)(a_1+a_4)/2} \nu^{-a_3} \bx^A
  \label{eq:Gl2-rhobasis}
\end{equation}
defines an isomorphic mapping between the $Gl_{p,q}(2)$~coalgebras whose
parameters are given in Eq.~(\ref{eq:gl2-pqchange}). By taking
Eqs.~(\ref{eq:Gl2-basisbis}) and~(\ref{eq:Gl2-h}) into account, and by
summing both
sides of Eq.~(\ref{eq:Gl2-rhobasis}) over~$a_2$ and~$a_3$, we also obtain the
relation
\begin{equation}
  \rho^{\nu}\left(\gamma^{a_1} e^{r(\alpha+\delta)} e^{s(\alpha-\delta)}
  \beta^{a_4}\right) = P^{(1-\nu)(a_1+a_4)/2} \gamma^{a_1} e^{r(\alpha+\delta)}
  e^{s(\alpha-\delta)/\nu} \beta^{a_4}, \label{eq:Gl2-rhoexp}
\end{equation}
valid for any $a_1$, $a_4 \in \N$, and any $r$,~$s \in \C$.\par
%
%
The action of~$\rho^{\nu}$, and of its inverse~$\rho_{\nu}$, on $\alpha$,
$\beta$,
$\gamma$,~$\delta$,
\begin{eqnarray}
  \rho^{\nu}(\alpha) & = & \case{\nu+1}{2\nu}\, \alpha +
\case{\nu-1}{2\nu}\, \delta,
         \qquad \rho^{\nu}(\beta) = P^{(1-\nu)/2} \beta, \nonumber \\
  \rho^{\nu}(\gamma) & = & P^{(1-\nu)/2} \gamma, \qquad \rho^{\nu}(\delta) =
         \case{\nu-1}{2\nu}\, \alpha + \case{\nu+1}{2\nu}\, \delta, \\
  \rho_{\nu}(\alpha) & = & \case{1+\nu}{2}\, \alpha + \case{1-\nu}{2}\, \delta,
         \qquad \rho_{\nu}(\beta) = P^{(\nu-1)/2} \beta, \nonumber \\
  \rho_{\nu}(\gamma) & = & P^{(\nu-1)/2} \gamma, \qquad \rho_{\nu}(\delta) =
         \case{1-\nu}{2}\, \alpha + \case{1+\nu}{2}\, \delta,
\end{eqnarray}
and on $a$, $b$, $c$, $d$, $\cal D$,
\begin{eqnarray}
  \rho^{\nu}(a) & = & {\cal D}^{(1-\nu)/(2\nu)} a, \qquad
        \rho^{\nu}(b) = \left(P^{\nu} {\cal D}\right)^{(1-\nu)/(2\nu)} b,
\nonumber \\
  \rho^{\nu}(c) & = & \left(P^{-\nu} {\cal D}\right)^{(1-\nu)/(2\nu)} c, \qquad
        \rho^{\nu}(d) = {\cal D}^{(1-\nu)/(2\nu)} d, \qquad
        \rho^{\nu}({\cal D}) = {\cal D}^{1/\nu}, \\
  \rho_{\nu}(a) & = & {\cal D}^{(\nu-1)/2} a, \qquad
        \rho_{\nu}(b) = (P {\cal D})^{(\nu-1)/2} b, \nonumber \\
  \rho_{\nu}(c) & = & \left(P^{-1} {\cal D}\right)^{(\nu-1)/2} c, \qquad
        \rho_{\nu}(d) = {\cal D}^{(\nu-1)/2} d, \qquad
        \rho_{\nu}({\cal D}) = {\cal D}^{\nu},    \label{eq:Gl2-rhoainv}
\end{eqnarray}
can be obtained as special cases of Eqs.~(\ref{eq:Gl2-rhobasis}),
and~(\ref{eq:Gl2-rhoexp}), respectively~\cite{footnote3}.\par
%
%
By using Eqs.~(\ref{eq:col-maps}), (\ref{eq:Gl2-alg}), (\ref{eq:Gl2-detcom}),
(\ref{eq:Gl2-S}), (\ref{eq:Gl2-gauss}), (\ref{eq:Gl2-alphadelta}),
(\ref{eq:Gl2-rhoexp}), and~(\ref{eq:Gl2-rhoainv}), it is straigthforward to
determine the coloured maps of~$Gl^c(2)$. The results are listed
in Appendix~A.\par
%
%
We shall now proceed to show that the relations obtained for
$\tm^{\nu}_{p,q,\lambda,\mu}$ can be rewritten in an alternative way. For such a
purpose, in analogy with Eq.~(\ref{eq:extT-rep}), let us define
\begin{equation}
  x(\lambda) \nudot y(\mu) \equiv \tm^{\nu}_{p,q,\lambda,\mu}(x \otimes y),
\end{equation}
where to avoid ambiguity, we specify the algebra to which each element
belongs by
the corresponding colour parameter, e.g., $x(\lambda) \in
Gl_{p^{(\lambda)},q^{(\lambda)}}(2)$.\par
%
%
The sixteen relations expressing $\tm^{\nu}_{p,q,\lambda,\mu}(x \otimes
y)$, $x$,
$y \in \{a, b, c, d\}$, as elements of~$Gl_{p^{(\nu)},q^{(\nu)}}(2)$, and
given in
Eqs.~(\ref{eq:Gl2-xy}), (\ref{eq:Gl2-txy}), can be combined into two
different sets.
The first one contains the relations
\begin{eqnarray}
  a(\lambda) \nudot b(\mu) & = & q^{(\lambda)}\, b(\mu) \nudot a(\lambda), \qquad
  a(\lambda) \nudot c(\mu) = p^{(\lambda)}\, c(\mu) \nudot a(\lambda),
\nonumber \\
  b(\lambda) \nudot d(\mu) & = & p^{(\mu)}\, d(\mu) \nudot a(\lambda), \qquad
  c(\lambda) \nudot d(\mu) = q^{(\mu)}\, d(\mu) \nudot c(\lambda), \nonumber \\
  b(\lambda) \nudot c(\mu) & = & \left(\frac{p^{(\lambda)}p^{(\mu)}}
          {q^{(\lambda)}q^{(\mu)}}\right)^{1/2}  c(\mu) \nudot b(\lambda),
\nonumber \\
  a(\lambda) \nudot d(\mu) - d(\mu) \nudot a(\lambda) & = & \left\{\left(
          q^{(\lambda)}q^{(\mu)}\right)^{1/2} - \left(p^{(\lambda)}p^{(\mu)}
          \right)^{-1/2}\right\}  b(\lambda) \nudot c(\mu),
\label{eq:Gl2-colalg1}
\end{eqnarray}
giving back the $Gl_{p,q}(2)$ defining relations~(\ref{eq:Gl2-alg}) in the
$\lambda$, $\mu$,~$\nu \to 1$ limit, while the second set is given by the
equations
\begin{eqnarray}
  a(\lambda) \nudot a(\mu) & = & a(\mu) \nudot a(\lambda), \qquad
         a(\lambda) \nudot b(\mu) = P^{(\mu-\lambda)/2}\, a(\mu) \nudot
b(\lambda),
         \nonumber \\
  a(\lambda) \nudot c(\mu) & = & P^{(\lambda-\mu)/2}\, a(\mu) \nudot c(\lambda),
         \qquad a(\lambda) \nudot d(\mu) = a(\mu) \nudot d(\lambda),
\nonumber \\
  b(\lambda) \nudot b(\mu) & = & P^{(\mu-\lambda)}\, b(\mu) \nudot b(\lambda),
         \qquad b(\lambda) \nudot c(\mu) = b(\mu) \nudot c(\lambda),
\nonumber \\
  b(\lambda) \nudot d(\mu) & = & P^{(\mu-\lambda)/2}\, b(\mu) \nudot d(\lambda),
         \qquad c(\lambda) \nudot c(\mu) = P^{\lambda-\mu}\, c(\mu) \nudot
         c(\lambda), \nonumber \\
  c(\lambda) \nudot d(\mu) & = & P^{(\lambda-\mu)/2}\, c(\mu) \nudot d(\lambda),
         \qquad d(\lambda) \nudot d(\mu) = d(\mu) \nudot d(\lambda),
         \label{eq:Gl2-colalg2}
\end{eqnarray}
which have no standard conterpart.\par
%
%
Similarly, the relations for the coloured multiplication involving the quantum
determinant, given in Eqs.~(\ref{eq:Gl2-Dx}), (\ref{eq:Gl2-xD}),
(\ref{eq:Gl2-DD}),
(\ref{eq:Gl2-tDx}), and~(\ref{eq:Gl2-txD}), can be rewritten as
\begin{eqnarray}
  {\cal D}(\lambda) \nudot a(\mu) & = & a(\mu) \nudot {\cal D}(\lambda), \qquad
          {\cal D}(\lambda) \nudot b(\mu) = P^{-2\lambda}\, b(\mu) \nudot
          {\cal D}(\lambda), \nonumber \\
  {\cal D}(\lambda) \nudot c(\mu) & = & P^{2\lambda}\, c(\mu) \nudot {\cal
          D}(\lambda), \qquad {\cal D}(\lambda) \nudot d(\mu) = d(\mu) \nudot
          {\cal D}(\lambda),
\end{eqnarray}
and
\begin{eqnarray}
  {\cal D}(\lambda) \nudot a(\mu) & = & ({\cal D}(\nu))^{(\lambda-\mu)/(2\nu)}
         \left({\cal D}(\mu) \nudot a(\lambda)\right), \nonumber \\
  {\cal D}(\lambda) \nudot b(\mu) & = & P^{(\mu-\lambda)/2}
         ({\cal D}(\nu))^{(\lambda-\mu)/(2\nu)} \left({\cal D}(\mu) \nudot
         b(\lambda)\right), \nonumber \\
  {\cal D}(\lambda) \nudot c(\mu) & = & P^{(\lambda-\mu)/2}
         ({\cal D}(\nu))^{(\lambda-\mu)/(2\nu)} \left({\cal D}(\mu) \nudot
         c(\lambda)\right), \nonumber \\
  {\cal D}(\lambda) \nudot d(\mu) & = & ({\cal D}(\nu))^{(\lambda-\mu)/(2\nu)}
         \left({\cal D}(\mu) \nudot d(\lambda)\right), \nonumber \\
  {\cal D}(\lambda) \nudot {\cal D}(\mu) & = & {\cal D}(\mu) \nudot
         {\cal D}(\lambda),
\end{eqnarray}
where the first set leads to Eq.~(\ref{eq:Gl2-detcom}) in the $\lambda$,
$\mu$,~$\nu \to 1$ limit, while the second has again no standard
counterpart.\par
%
%
The simplicity of Eqs.~(\ref{eq:Gl2-colalg1}), and~(\ref{eq:Gl2-colalg2})
should be
stressed. It is also remarkable that these equations are entirely
independent of the
final Hopf algebra colour parameter~$\nu$. Such a property reflects the
fact that
the two sets of equations (\ref{eq:Gl2-colalg1}),~(\ref{eq:Gl2-colalg2}) are
equivalent to the single matrix equation~(\ref{eq:colRTT-defrep}), where
$R^{\lambda,\mu}_{p,q}$ is given by Eq.~(\ref{eq:gl2-colR}), and
$T^{\lambda}_{p,q}
= T_{p^{(\lambda)},q^{(\lambda)}} = \left(\begin{array}{cc} a(\lambda) & b(\lambda)
\\ c(\lambda) & d(\lambda) \end{array}\right)$. It can indeed be checked by
expressing Eq.~(\ref{eq:colRTT-defrep}) in elementwise form that the resulting
independent relations can be grouped together in two sets, given by
Eqs.~(\ref{eq:Gl2-colalg1}), and~(\ref{eq:Gl2-colalg2}). Since, moreover, the
coalgebraic structure of~$Gl^c(2)$ coincides with that of each individual
standard
quantum group~$Gl_{p,q}(2)$, for $(p,q)$ running over~$\cal Q$, we did show in a
very explicit way that our coloured $Gl^c(2)$ is equivalent to the coloured
extension of~$Gl_{p,q}(2)$, defined through coloured $RTT$-relations,
according to
Basu-Mallick's prescription~\cite{basu}.\par
%
%
By the way, it is worth noting that since parameter~$Q$ is left unchanged by the
colour group transformations, the elements~${\cal T}^{\lambda}_{p,q}$ of the
coloured universal $\cal T$-matrix, considered in
Definition~\ref{def-colT}, have
no explicit $\lambda$-dependence, so that they are all given by
Eq.~(\ref{eq:Gl2-univTbis}) with $\hJ_0$, $\hJ_{\pm}$, $\hZ \in
U_{p^{(\lambda)},
q^{(\lambda)}}(gl(2))$, and $h$, $\th$, $\beta$, $\gamma \in Gl_{p^{(\lambda)},
q^{(\lambda)}}(2)$.\par
%
%
As a final point of this Section, let us mention that the coloured maps
$\tm^{\nu}_{p,q,\lambda,\mu}$, $\tiota^{\nu}_{p,q}$,
and~$\tS^{\nu}_{p,q,\mu}$ may
be extended from~$Gl^c(2)$ to the whole dual of~$U^c(gl(2))$ by considering the
enveloping algebras $U_{p,q}(\{\alpha,\beta,\gamma,\delta\})$ for $(p,q)$
running
over~$\cal Q$. For the coloured multiplication, for instance, one finds that the
sixteen relations giving $\tm^{\nu}_{p,q,\lambda,\mu}(x \otimes y)$, $x$, $y \in
\{\alpha,\beta,\gamma,\delta\}$, can be combined into six generalized
commutation relations
\begin{eqnarray}
  [\alpha(\lambda), \beta(\mu)]^{\nu} & = & P^{(\mu-\nu)/2} (\varphi - \lambda
         \theta) \beta(\nu), \qquad  [\alpha(\lambda), \gamma(\mu)]^{\nu} =
         P^{(\mu-\nu)/2} (\varphi + \lambda \theta) \gamma(\nu), \nonumber
\\[0cm]
  [\delta(\lambda), \beta(\mu)]^{\nu} & = & P^{(\mu-\nu)/2} (\varphi + \lambda
         \theta) \beta(\nu), \qquad  [\delta(\lambda), \gamma(\mu)]^{\nu} =
         P^{(\mu-\nu)/2} (\varphi - \lambda \theta) \gamma(\nu), \nonumber
\\[0cm]
  [\alpha(\lambda), \delta(\mu)]^{\nu} & = & 0, \qquad  [\beta(\lambda),
         \gamma(\mu)]^{\nu} = 0,   \label{eq:Gl2-colLie1}
\end{eqnarray}
where
\begin{equation}
  [x(\lambda), y(\mu)]^{\nu} \equiv x(\lambda) \nudot y(\mu) - y(\mu) \nudot
  x(\lambda),
\end{equation}
and ten additional relations
\begin{eqnarray}
  [\alpha(\lambda), \alpha(\mu)]^{\nu} & = & [\beta(\lambda),
\beta(\mu)]^{\nu} =
           [\gamma(\lambda), \gamma(\mu)]^{\nu} = [\delta(\lambda),
           \delta(\mu)]^{\nu} = 0, \nonumber \\
  \alpha(\lambda) \nudot \beta(\mu) & = & P^{(\mu-\lambda)/2} (2\mu)^{-1}
           \left((\lambda+\mu) \alpha(\mu) \nudot \beta(\lambda) - (\lambda-\mu)
           \delta(\mu) \nudot \beta(\lambda)\right), \nonumber \\
  \gamma(\lambda) \nudot \alpha(\mu) & = & P^{(\lambda-\mu)/2} (2\lambda)^{-1}
           \left((\lambda+\mu) \gamma(\mu) \nudot \alpha(\lambda) +
(\lambda-\mu)
           \gamma(\mu) \nudot \delta(\lambda)\right), \nonumber \\
  \alpha(\lambda) \nudot \delta(\mu) & = & \alpha(\mu) \nudot \delta(\lambda)
           + (\lambda-\mu) (\lambda+\mu)^{-1} \left(\alpha(\mu) \nudot
           \alpha(\lambda) - \delta(\mu) \nudot \delta(\lambda)\right),
\nonumber \\
  \gamma(\lambda) \nudot \beta(\mu) & = & \gamma(\mu) \nudot \beta(\lambda),
           \nonumber \\
  \delta(\lambda) \nudot \beta(\mu) & = & P^{(\mu-\lambda)/2} (2\mu)^{-1}
           \left(- (\lambda-\mu) \alpha(\mu) \nudot \beta(\lambda) +
(\lambda+\mu)
           \delta(\mu) \nudot \beta(\lambda)\right), \nonumber \\
  \gamma(\lambda) \nudot \delta(\mu) & = & P^{(\lambda-\mu)/2} (2\lambda)^{-1}
           \left((\lambda-\mu) \gamma(\mu) \nudot \alpha(\lambda) +
(\lambda+\mu)
           \gamma(\mu) \nudot \delta(\lambda)\right).   \label{eq:Gl2-colLie2}
\end{eqnarray}
Whenever $\lambda$, $\mu$, $\nu \to 1$, the former set reproduces the defining
relations~(\ref{eq:Gl2-Lie}) of~$U_{p,q}(\{\alpha,\beta,\gamma,\delta\})$,
whereas the latter has no standard counterpart. Let us stress that contrary to
what happens in Eqs.~(\ref{eq:Gl2-colalg1}) and~(\ref{eq:Gl2-colalg2}), there
appears an explicit $\nu$-dependence in Eqs.~(\ref{eq:Gl2-colLie1}),
and~(\ref{eq:Gl2-colLie2}).\par
%
%
\section{THE COLOURED TWO-PARAMETER QUEA $U^c(gl(1/1))$ AND ITS DUAL
$Gl^c(1/1)$}
\label{sec:gl(1/1)}
\setcounter{equation}{0}
One of the simplest examples of dual pairs made of a QUEA of a Lie superalgebra
and of a quantum supergroup consists in the two-parameter deformations
of~$U(gl(1/1))$ and $Gl(1/1)$~\cite{dabrowski}--\cite{chakra95}, which find
some interesting applications to the multivariable Alexander-Conway polynomial
and the free fermion model~\cite{kauffman}, and, in another context, to the
XY~quantum chain in a magnetic field~\cite{hinrichsen}.\par
%
%
In~II, a coloured QUEA was constructed by starting from the two-parameter
deformation of~$U(gl(1/1))$ in the Bednar {\it et al\/}
formulation~\cite{bednar,
burdik92b}. Moreover, Chakrabarti and Jagannathan~\cite{chakra96b} recently
extended the Fronsdal and Galindo universal {\cal T}-matrix
formalism~\cite{fronsdal93} to the pair of Hopf superalgebras~$U_{p,q}(gl(1/1))$
and~$Gl_{p,q}(1/1)$. So all the ingredients needed to build an example of
coloured
Hopf dual in the graded case are available. To solve such a problem, we
shall use
the approach to~$U_{p,q}(gl(1/1))$ employed in Ref.~\cite{chakra96b}, instead of
that of Refs.~\cite{bednar,burdik92b}; hence, the results
for~$U^c(gl(1/1))$ derived
here will slightly differ from those obtained in~II.\par
%
%
\subsection{The dual Hopf superalgebras $U_{p,q}(gl(1/1))$ and $Gl_{p,q}(1/1)$}
The quantum superalgebra $U_{p,q}(gl(1/1))$ is defined~\cite{chakra96b} as the
algebra generated by two even elements $\{J_0, Z\}$, and two odd ones $\{J_+,
J_-\}$, subject to the relations
\begin{equation}
  \left[J_0, J_{\pm}\right] = \pm J_{\pm}, \qquad \left\{J_+, J_-\right\} =
  [2Z]_Q, \qquad J_{\pm}^2 = 0, \qquad \left[Z, J_0\right] = \left[Z,
J_{\pm}\right] =
  0,
\end{equation}
where $[X]_t$, and $Q$,~$P$ are still defined by Eqs.~(\ref{eq:box})
and~(\ref{eq:QP}), respectively. As in the $U_{p,q}(gl(2))$ case, the
coalgebraic
structure depends upon both $P$ and~$Q$.\par
%
%
The duality relationship between $U_{p,q}(gl(1/1))$ and~$Gl_{p,q}(1/1)$ is most
easily described in terms of another set of $U_{p,q}(gl(1/1))$ generators
$\{\hJ_0,
\hJ_{\pm}, \hZ\}$, with $\hJ_0$,~$\hZ$ even, and $\hJ_{\pm}$ odd. The new
generators are defined in terms of the old ones by
\begin{equation}
  \hJ_0 = J_0, \qquad \hJ_{\pm} = J_{\pm}\, Q^{\pm Z}\, P^{Z - \half},
  \qquad \hZ = Z,
\end{equation}
and satisfy the superalgebra (anti)commutation relations
\begin{eqnarray}
  \left[\hJ_0, \hJ_{\pm}\right] & = & \pm \hJ_{\pm}, \qquad \left\{\hJ_+,
         \hJ_-\right\} = \frac{p^{2\hZ}-q^{-2\hZ}}{p-q^{-1}}, \qquad
\hJ_{\pm}^2 = 0,
         \nonumber \\
  \left[\hZ, \hJ_0\right] & = & \left[\hZ, \hJ_{\pm}\right] = 0.
  \label{eq:gl1/1-algbis}
\end{eqnarray}
The remaining Hopf superalgebra maps are given by
\begin{eqnarray}
  \Delta_{p,q}\left(\hJ_0\right) & = & \hJ_0 \otimes 1 + 1 \otimes \hJ_0, \qquad
         \Delta_{p,q}\left(\hZ\right) = \hZ \otimes 1 + 1 \otimes \hZ,
\nonumber \\
  \Delta_{p,q}\left(\hJ_+\right) & = & \hJ_+ \otimes p^{2\hZ} + 1 \otimes \hJ_+,
         \qquad \Delta_{p,q}\left(\hJ_-\right) = \hJ_- \otimes 1
         + q^{-2\hZ} \otimes \hJ_-, \nonumber \\
  \epsilon_{p,q}(X) & = & 0, \qquad X \in \{\hJ_0, \hJ_{\pm}, \hZ\}, \nonumber \\
  S_{p,q}\left(\hJ_0\right) & = & - \hJ_0, \qquad S_{p,q}\left(\hZ\right) =
- \hZ,
         \nonumber \\
  S_{p,q}\left(\hJ_+\right) & = & - p^{-2\hZ} \hJ_+, \qquad
         S_{p,q}\left(\hJ_-\right) = - q^{2\hZ} \hJ_-,
\end{eqnarray}
and the universal $\cal R$-matrix can be written as~\cite{burdik92b,chakra95}
\begin{equation}
  {\cal R}_{p,q} = p^{2\hZ \otimes \hJ_0} \, q^{2\hJ_0 \otimes \hZ} \left\{ 1
  \otimes 1 - \left(p - q^{-1}\right) \hJ_+ \otimes \hJ_-\right\}.
  \label{eq:gl1/1-univR}
\end{equation}
\par
%
%
In the $2 \times 2$ defining representation~$D_{p,q}$ of $U_{p,q}(gl(1/1))$,
given by Eq.~(\ref{eq:gl2-defrep}), the universal $\cal
R$-matrix~(\ref{eq:gl1/1-univR}) is represented by the $4 \times 4$~matrix
\begin{equation}
  R_{p,q} \equiv \left(D_{p,q} \otimes D_{p,q}\right) \left({\cal
R}_{p,q}\right) =
              \left(\begin{array}{cccc}
                      Q & 0        & 0              & 0 \\[0.1cm]
                      0 & P^{-1} & Q - Q^{-1} & 0 \\[0.1cm]
                      0 & 0        & P              & 0 \\[0.1cm]
                      0 & 0        & 0              & Q^{-1}
              \end{array} \right).   \label{eq:gl1/1-R}
\end{equation}
In deriving Eq.~(\ref{eq:gl1/1-R}), use is made of
Eq.~(\ref{eq:tensorprod}) with
the $\Z_2$-grade of the first (resp.~second) row or column of~$D_{p,q}(X)$
defined
as zero (resp.~one).\par
%
%
The defining $T$-matrix of the corresponding quantum supergroup~$Gl_{p,q}(1/1)$
is still given by Eq.~(\ref{eq:Gl2-T}), where $a$,~$d$, and $b$,~$c$ are
now even
and odd, respectively, and satisfy the relations
\begin{eqnarray}
  ab & = & p^{-1} ba, \qquad ac = q^{-1} ca, \qquad bd = pdb, \qquad cd = qdc,
          \nonumber \\
  bc & = & - (p/q) cb, \qquad ad - da = \left(q - p^{-1}\right) bc, \qquad
b^2 = c^2 =
          0.     \label{eq:Gl1/1-alg}
\end{eqnarray}
The latter follow from the $RTT$-relations corresponding to the
$R$-matrix~(\ref{eq:gl1/1-R}), with convention~(\ref{eq:A1-A2}) taken into
account. The coalgebra maps are given by Eq.~(\ref{eq:Gl2-coalg}), while the
antipode one is
\begin{equation}
  \tS_{p,q}(T_{p,q}) = (T_{p,q})^{-1} =
  \left(\begin{array}{cc}
            a^{-1} + a^{-1}bd^{-1}ca^{-1} & - a^{-1}bd^{-1} \\[0.1cm]
            - d^{-1}ca^{-1}                      & d^{-1} + d^{-1}ca^{-1}bd^{-1}
           \end{array}\right),
\end{equation}
where $a$ and $d$ are assumed invertible. The quantum superdeterminant, defined
by
\begin{equation}
  {\cal D} \equiv ad^{-1} - bd^{-1}cd^{-1},
\end{equation}
is both central and group-like, with $\tS_{p,q}({\cal D}) = {\cal D}^{-1} =
da^{-1} +
ba^{-1}ca^{-1}$.\par
%
%
Gauss decomposition~(\ref{eq:Gl2-gauss}) still holds, but now $a$ and $\hd$ are
invertible even elements, while $\beta$ and~$\gamma$ are odd. The antipode map
and the quantum superdeterminant may be rewritten as
\begin{equation}
  \tS_{p,q}(T_{p,q}) =
  \left(\begin{array}{cc}
            a^{-1} - p^{-1}qca^{-2}\hd^{-1}b & - p^{-1}a^{-1}\hd^{-1}b \\[0.1cm]
            - qca^{-1}\hd^{-1}                       & \hd^{-1}
           \end{array}\right),    \label{eq:Gl1/1-S}
\end{equation}
and
\begin{equation}
  {\cal D} = a \hd^{-1},
\end{equation}
respectively.\par
%
%
By assuming that the superalgebra can be augmented with the logarithms of~$a$
and~$\hd$, and by using the maps
\begin{equation}
  a = e^{\alpha}, \qquad \hd = e^{\delta},   \label{eq:Gl1/1-alphadelta}
\end{equation}
and the reparametrization
\begin{equation}
  p = e^{-\omega}, \qquad q = e^{-\eta},
\end{equation}
the quantum superdeterminant is transformed into ${\cal D} = \exp(\alpha -
\delta)$, and the new variables $\{\alpha, \beta, \gamma, \delta\}$ satisfy the
(anti)commutation relations of a solvable Lie superalgebra
\begin{eqnarray}
  [\alpha, \beta] & = & \omega \beta, \qquad [\alpha, \gamma] = \eta
\gamma, \qquad
           [\delta, \beta] = \omega \beta, \qquad [\delta, \gamma] = \eta
\gamma,
           \nonumber \\[0cm]
  [\alpha, \delta] & = & 0, \qquad \{\beta, \gamma\} = 0, \qquad \beta^2 =
0, \qquad
           \gamma^2 = 0,    \label{eq:Gl1/1-Lie}
\end{eqnarray}
with a noncocommutative coproduct structure~\cite{chakra96b}. This shows that
$Gl_{p,q}(1/1)$ can be embedded into the enveloping algebra $U_{p,q}(\{\alpha,
\beta,\gamma,\delta\})$ of such a Lie superalgebra. According to Chakrabarti and
Jagannathan~\cite{chakra96b}, $U_{p,q}(\{\alpha, \beta, \gamma, \delta\})$
is dual
to the QUEA~$U_{p,q}(gl(1/1))$~\cite{footnote4}.\par
%
%
Dual bases of~$U_{p,q}(gl(1/1))$ and~$Gl_{p,q}(1/1)$ are given
by~\cite{chakra96b}
\begin{equation}
  X_A = \hJ_-^{a_1} \frac{H^{a_2}}{a_2!} \frac{\tH^{a_3}}{a_3!}
\hJ_+^{a_4}, \qquad
  x^A = \gamma^{a_1} \alpha^{a_2} \delta^{a_3} \beta^{a_4},
  \label{eq:gl1/1-basis}
\end{equation}
respectively, where $H \equiv \hZ + \hJ_0$, $\tH \equiv \hZ - \hJ_0$,
$A = (a_1, a_2, a_3, a_4)$, $a_1$, $a_4 \in \{0,1\}$, and $a_2$,~$a_3 \in
\N$. From
Eqs.~(\ref{eq:T}), and~(\ref{eq:gl1/1-basis}), the universal $\cal T$-matrix
of~$Gl_{p,q}(1/1)$ can be written as
\begin{equation}
  {\cal T}_{p,q} = \exp\left(\gamma \hJ_-\right) \exp\left(\alpha H +
  \delta \tH\right) \exp\left(\beta \hJ_+\right),   \label{eq:Gl1/1-univT}
\end{equation}
in terms of standard exponentials.\par
%
%
In the next Subsection, it will prove convenient to consider another set of dual
bases
\begin{eqnarray}
  \bX_A & = & \sum_B c_A^B X_B = \hJ_-^{a_1} \frac{\hJ_0^{a_2}}{a_2!}
           \frac{\hZ^{a_3}}{a_3!} \hJ_+^{a_4}, \label{eq:gl1/1-basisbis} \\
  \bx^A & = & \sum_B \left(c^{-1}\right)^A_B x^B = \gamma^{a_1} h^{a_2}
\th^{a_3}
           \beta^{a_4},
\end{eqnarray}
where $A$ has the same meaning as in Eq.~(\ref{eq:gl1/1-basis}),
\begin{equation}
  h \equiv \alpha - \delta, \qquad \th \equiv \alpha + \delta,
\end{equation}
and
\begin{eqnarray}
  c_A^B & = & 2^{-a_2-a_3} (-1)^{a_2-b_2} \left(c^{-1}\right)_A^B \nonumber \\
  & = & \delta_{a_1}^{b_1}\, \delta_{a_2+a_3}^{b_2+b_3}\, \delta_{a_4}^{b_4}\,
         2^{-a_2-a_3} \sum_{t=\max(0,b_2-a_2)}^{\min(b_2,a_3)} (-1)^{a_2-b_2+t}
         \bin{b_2}{t} \bin{b_3}{a_3-t}.
\end{eqnarray}
In terms of them, the universal $\cal T$-matrix~(\ref{eq:Gl1/1-univT}) can be
rewritten as
\begin{equation}
  {\cal T}_{p,q} = \exp\left(\gamma \hJ_-\right) \exp\left(h \hJ_0 + \th
\hZ\right)
  \exp\left(\beta \hJ_+\right).
\end{equation}
\par
%
%
\subsection{The coloured Hopf superalgebra $U^c(gl(1/1))$ and its dual
$Gl^c(1/1)$}
The defining relations~(\ref{eq:gl1/1-algbis}) of the
$U_{p,q}(gl(1/1))$~superalgebra are left invariant under the grade-preserving
transformations
\begin{equation}
  \sigma^{\nu}\left(\hJ_0\right) = \hJ_0, \qquad
  \sigma^{\nu}\left(\hJ_{\pm}\right) = \left(\frac{p^{\nu}-q^{-\nu}}{p-q^{-1}}
  \right)^{1/2} \hJ_{\pm}, \qquad \sigma^{\nu}\left(\hZ\right) = \nu \hZ,
  \label{eq:gl1/1-sigma}
\end{equation}
where $\nu \in {\cal C} = \C \setminus \{0\}$, provided $p$ and~$q$ are
replaced by
their $\nu$th~powers, $p^{\nu}$ and~$q^{\nu}$ (and the same for $P$
and~$Q$). The
$\sigma^{\nu}$'s are therefore isomorphic mappings between the
$U_{p,q}(gl(1/1))$ and $U_{p^{\nu},q^{\nu}}(gl(1/1))$~superalgebras, and define
a colour group ${\cal G} = Gl(1,\C)$.\par
%
%
By choosing the parameter set ${\cal Q} = (\C \setminus \{0\}) \times (\C
\setminus \{0\})$, we obtain a coloured quasitriangular Hopf
superalgebra~$U^c(gl(1/1))$, whose coloured comultiplication, counit,
antipode, and
$\cal R$-matrix are given by
\begin{eqnarray}
  \Delta^{\lambda,\mu}_{p,q,\nu}\left(\hJ_0\right) & = & \hJ_0 \otimes 1 + 1
           \otimes \hJ_0, \qquad
\Delta^{\lambda,\mu}_{p,q,\nu}\left(\hZ\right) =
           \frac{\lambda}{\nu}\, \hZ \otimes 1 + \frac{\mu}{\nu}\, 1
\otimes \hZ,
           \nonumber \\
  \Delta^{\lambda,\mu}_{p,q,\nu}\left(\hJ_+\right) & = & A^{\lambda}_{\nu}\,
           \hJ_+ \otimes p^{2\mu\hZ} + A^{\mu}_{\nu}\, 1 \otimes \hJ_+, \qquad
           \Delta^{\lambda,\mu}_{p,q,\nu}\left(\hJ_-\right) =
A^{\lambda}_{\nu}\,
           \hJ_- \otimes 1 + A^{\mu}_{\nu}\, q^{-2\lambda\hZ} \otimes \hJ_-,
           \nonumber \\
  \epsilon_{p,q,\nu}(X) & = & 0, \qquad X \in \{\hJ_0, \hJ_{\pm}, \hZ\},
\nonumber \\
  S^{\mu}_{p,q,\nu}\left(\hJ_0\right) & = & - \hJ_0, \qquad
           S^{\mu}_{p,q,\nu}\left(\hJ_+\right) = - A^{\mu}_{\nu}\,
p^{-2\mu\hZ} \hJ_+,
           \nonumber \\
  S^{\mu}_{p,q,\nu}\left(\hJ_-\right) & = & - A^{\mu}_{\nu}\, q^{2\mu\hZ} \hJ_-,
           \qquad S^{\mu}_{p,q,\nu}\left(\hZ\right) = - \frac{\mu}{\nu}
\hZ, \nonumber
           \\
  {\cal R}^{\lambda,\mu}_{p,q} & = & p^{2\lambda \hZ \otimes \hJ_0}\,
           q^{2\mu \hJ_0 \otimes \hZ} \left\{1 \otimes 1 -
\sqrt{\left(p^{\lambda} -
           q^{-\lambda}\right) \left(p^{\mu}-q^{-\mu}\right)}\, \hJ_+
\otimes \hJ_-
           \right\},
\end{eqnarray}
where $A^{\lambda}_{\nu} \equiv \bigl((p^{\lambda}-q^{-\lambda}) \big/
(p^{\nu}-q^{-\nu})\bigr)^{1/2}$.\par
%
%
In the defining representation~(\ref{eq:gl2-defrep}), the coloured
universal $\cal
R$-matrix is represented by the $4 \times 4$~matrix
\begin{equation}
  R^{\lambda,\mu}_{p,q} =
     \left(\begin{array}{cccc}
          \ss Q^{(\lambda+\mu)/2} P^{(\lambda-\mu)/2} & \ss 0 & \ss 0 & \ss 0
               \\[0.1cm]
          \ss 0 & \ss Q^{-(\lambda-\mu)/2} P^{-(\lambda+\mu)/2} & \ss
               Q^{-(\lambda-\mu)/2} B^{\lambda,\mu}  & \ss 0 \\[0.1cm]
          \ss 0 & \ss 0 & \ss Q^{(\lambda-\mu)/2} P^{(\lambda+\mu)/2} & \ss 0
               \\[0.1cm]
          \ss 0 & \ss 0 & \ss 0 & \ss Q^{-(\lambda+\mu)/2} P^{-(\lambda-\mu)/2}
     \end{array} \right),   \label{eq:gl1/1-colR}
\end{equation}
with $B^{\lambda,\mu} \equiv \sqrt{\left(Q^{\lambda} - Q^{-\lambda}\right)
\left(Q^{\mu}-Q^{-\mu}\right)}$. Such a matrix is a solution of the coloured
graded YBE, and gives back matrix~(\ref{eq:gl1/1-R}) whenever
$\lambda$,~$\mu\to 1$.\par
%
%
Under transformation~(\ref{eq:gl1/1-sigma}), the basis elements~$\bX_A$
of~$U_{p,q}(gl(1/1))$, defined in Eq.~(\ref{eq:gl1/1-basisbis}), become
\begin{equation}
  \sigma^{\nu}\left(\bX_A\right) = \left(A^{\nu}_1\right)^{(a_1+a_4)/2}
\nu^{a_3}
  \bX_A.
\end{equation}
Hence, from Eqs.~(\ref{eq:rho}) and~(\ref{eq:pairing}), it follows that
\begin{equation}
  \rho^{\nu}\left(\bx^A\right) = \left(A^1_{\nu}\right)^{(a_1+a_4)/2} \nu^{-a_3}
  \bx^A,   \label{eq:Gl1/1-rhobasis}
\end{equation}
and
\begin{equation}
  \rho^{\nu}\left(\gamma^{a_1} e^{r(\alpha-\delta)} e^{s(\alpha+\delta)}
  \beta^{a_4}\right) = \left(A^1_{\nu}\right)^{(a_1+a_4)/2} \gamma^{a_1}
  e^{r(\alpha-\delta)} e^{s(\alpha+\delta)/\nu} \beta^{a_4},
\label{eq:Gl1/1-rhoexp}
\end{equation}
where $r$,~$s \in \C$.\par
%
%
As special cases of Eqs.~(\ref{eq:Gl1/1-rhobasis})
and~(\ref{eq:Gl1/1-rhoexp}), we
obtain~\cite{footnote3}
\begin{eqnarray}
  \rho^{\nu}(\alpha) & = & \case{1+\nu}{2\nu}\, \alpha + \case{1-\nu}{2\nu}\,
         \delta, \qquad \rho^{\nu}(\beta) = A^1_{\nu}\, \beta, \nonumber \\
  \rho^{\nu}(\gamma) & = & A^1_{\nu}\, \gamma, \qquad \rho^{\nu}(\delta) =
         \case{1-\nu}{2\nu}\, \alpha + \case{1+\nu}{2\nu}\, \delta, \\
  \rho_{\nu}(\alpha) & = & \case{\nu+1}{2}\, \alpha + \case{\nu-1}{2}\, \delta,
         \qquad \rho_{\nu}(\beta) = A^{\nu}_1\, \beta, \nonumber \\
  \rho_{\nu}(\gamma) & = & A^{\nu}_1\, \gamma, \qquad \rho_{\nu}(\delta) =
         \case{\nu-1}{2}\, \alpha + \case{\nu+1}{2}\, \delta,
\end{eqnarray}
and
\begin{eqnarray}
  \rho^{\nu}(a) & = & {\cal D}^{(\nu-1)/(2\nu)} a^{1/\nu}, \qquad
        \rho^{\nu}(b) = A^1_{\nu}\, {\cal D}^{(\nu-1)/(2\nu)} a^{(1-\nu)/\nu}b,
        \nonumber \\
  \rho^{\nu}(c) & = & q^{1-\nu} A^1_{\nu}\, {\cal D}^{(\nu-1)/(2\nu)}
a^{(1-\nu)/\nu}
        c, \nonumber \\
  \rho^{\nu}(d) & = & {\cal D}^{(\nu-1)/(2\nu)} a^{(1-\nu)/\nu} \left\{d - p
        \left(A^{\nu-1}_{\nu}\right)^2 c a^{-1} b\right\}, \qquad
        \rho^{\nu}({\cal D}) = {\cal D}, \\
  \rho_{\nu}(a) & = & {\cal D}^{(1-\nu)/2} a^{\nu}, \qquad
        \rho_{\nu}(b) = A^{\nu}_1\, {\cal D}^{(1-\nu)/2} a^{\nu-1} b,
\nonumber \\
  \rho_{\nu}(c) & = & q^{\nu-1} A^{\nu}_1\, {\cal D}^{(1-\nu)/2} a^{\nu-1} c,
        \nonumber \\
  \rho_{\nu}(d) & = & {\cal D}^{(1-\nu)/2} a^{\nu-1} \left\{d - p^{\nu}
        \left(A^{1-\nu}_1\right)^2 c a^{-1} b\right\}, \qquad
        \rho_{\nu}({\cal D}) = {\cal D}.    \label{eq:Gl1/1-rhoainv}
\end{eqnarray}
\par
%
%
By proceeding as in Sec.~\ref{sec:gl(2)}, it is straightforward to determine the
coloured maps of~$Gl^c(1/1)$, listed in Appendix~B, and to combine the
results for
the coloured multiplication $\tm^{\nu}_{p,q,\lambda,\mu}(x \otimes y)$,
$x$, $y \in
\{a, b, c, d\}$, given in Eqs.~(\ref{eq:Gl1/1-xy})--(\ref{eq:Gl1/1-dd}),
and~(\ref{eq:Gl1/1-Cxy}), into two sets of eight relations each,
\begin{eqnarray}
  a(\lambda) \nudot b(\mu) & = & p^{-\lambda} b(\mu) \nudot a(\lambda), \qquad
          a(\lambda) \nudot c(\mu) = q^{-\lambda} c(\mu) \nudot a(\lambda),
          \nonumber \\
  b(\lambda) \nudot d(\mu) & = & p^{\mu}\, d(\mu) \nudot a(\lambda), \qquad
          c(\lambda) \nudot d(\mu) = q^{\mu}\, d(\mu) \nudot c(\lambda),
\nonumber \\
  b(\lambda) \nudot c(\mu) & = & - p^{\mu} q^{-\lambda}  c(\mu) \nudot
b(\lambda),
          \qquad b(\lambda) \nudot b(\mu) = c(\lambda) \nudot c(\mu) = 0,
\nonumber
          \\
  a(\lambda) \nudot d(\mu) - d(\mu) \nudot a(\lambda) & = & p^{-\mu} q^{\lambda}
          \sqrt{\left(p^{\lambda}-q^{-\lambda}\right)
          \left(p^{\mu}-q^{-\mu}\right)}\, b(\lambda) \nudot c(\mu),
          \label{eq:Gl1/1-colalg1}
\end{eqnarray}
and
\begin{eqnarray}
  a(\lambda) \nudot a(\mu) & = & a(\mu) \nudot a(\lambda), \qquad
         a(\lambda) \nudot b(\mu) = A^{\mu}_{\lambda}\, a(\mu) \nudot
b(\lambda),
         \nonumber \\
  a(\lambda) \nudot c(\mu) & = & q^{\mu-\lambda} A^{\mu}_{\lambda}\, a(\mu)
         \nudot c(\lambda), \nonumber \\
  a(\lambda) \nudot d(\mu) - a(\mu) \nudot d(\lambda) & = & q^{\lambda-\mu}
         A^{\lambda-\mu}_{\lambda} A^{\lambda-\mu}_{\mu}\, b(\lambda) \nudot
         c(\mu) , \nonumber \\
  b(\lambda) \nudot c(\mu) & = & (pq)^{\mu-\lambda} b(\mu) \nudot c(\lambda),
         \qquad b(\lambda) \nudot d(\mu) = p^{\mu-\lambda} A^{\lambda}_{\mu}\,
         b(\mu) \nudot d(\lambda), \nonumber \\
  c(\lambda) \nudot d(\mu) & = & A^{\lambda}_{\mu}\, c(\mu) \nudot d(\lambda),
         \qquad d(\lambda) \nudot d(\mu) = d(\mu) \nudot d(\lambda).
         \label{eq:Gl1/1-colalg2}
\end{eqnarray}
\par
%
%
Similarly, the relations for the coloured multiplication involving the quantum
superdeterminant, given in Eqs.~(\ref{eq:Gl1/1-Dx})--(\ref{eq:Gl1/1-DD}),
and~(\ref{eq:Gl1/1-CDx}), can be rewritten as
\begin{equation}
  {\cal D}(\lambda) \nudot x(\mu) = x(\mu) \nudot {\cal D}(\lambda), \qquad
  x \in \{a, b, c ,d\},   \label{eq:Gl1/1-coldet1}
\end{equation}
and
\begin{eqnarray}
  {\cal D}(\lambda) \nudot a(\mu) & = & ({\cal D}(\nu))^{(\lambda-\mu)/(2\nu)}
         (a(\nu))^{(\mu-\lambda)/\nu} \left({\cal D}(\mu) \nudot
a(\lambda)\right),
         \nonumber \\
  {\cal D}(\lambda) \nudot b(\mu) & = & A^{\mu}_{\lambda}\,
         ({\cal D}(\nu))^{(\lambda-\mu)/(2\nu)} (a(\nu))^{(\mu-\lambda)/\nu}
         \left({\cal D}(\mu) \nudot b(\lambda)\right), \nonumber \\
  {\cal D}(\lambda) \nudot c(\mu) & = & q^{\mu-\lambda} A^{\mu}_{\lambda}\,
         ({\cal D}(\nu))^{(\lambda-\mu)/(2\nu)} (a(\nu))^{(\mu-\lambda)/\nu}
         \left({\cal D}(\mu) \nudot c(\lambda)\right), \nonumber \\
  {\cal D}(\lambda) \nudot d(\mu) & = & q^{\mu-\lambda}
         \left(A^{\mu-\nu}_{\lambda-\nu}\right)^2
         ({\cal D}(\nu))^{(\lambda-\mu)/(2\nu)} (a(\nu))^{(\mu-\lambda)/\nu}
         \left({\cal D}(\mu) \nudot d(\lambda)\right) \nonumber \\
  & & \mbox{} + p^{\mu-\nu} \left(A^{\lambda-\mu}_{\lambda-\nu}\right)^2
         ({\cal D}(\nu))^{(3\nu-\mu)/(2\nu)} (a(\nu))^{(\mu-\lambda)/\nu}
d(\nu),
         \nonumber \\
  {\cal D}(\lambda) \nudot {\cal D}(\mu) & = & {\cal D}(\mu) \nudot
         {\cal D}(\lambda).   \label{eq:Gl1/1-coldet2}
\end{eqnarray}
\par
%
%
In the $\lambda$, $\mu$, $\nu \to 1$~limit, Eqs.~(\ref{eq:Gl1/1-colalg1}),
and~(\ref{eq:Gl1/1-coldet1}) give back the $Gl_{p,q}(1/1)$ defining
relations~(\ref{eq:Gl1/1-alg}) and the central property of $\cal D$, while
Eqs.~(\ref{eq:Gl1/1-colalg2}), and~(\ref{eq:Gl1/1-coldet2}) have no standard
counterpart.\par
%
%
The two $\nu$-independent sets of equations, contained
in~(\ref{eq:Gl1/1-colalg1}) and~(\ref{eq:Gl1/1-colalg2}), can be explicitly
shown
to be equivalent to the single matrix equation~(\ref{eq:colRTT-defrep}), where
$R^{\lambda,\mu}_{p,q}$ is given by Eq.~(\ref{eq:gl1/1-colR}). Moreover, as in
Sec.~\ref{sec:gl(2)}, the elements~${\cal T}^{\lambda}_{p,q}$ of the coloured
universal $\cal T$-matrix have no explicit $\lambda$-dependence.\par
%
%
The counterparts of Eqs.~(\ref{eq:Gl2-colLie1}) and~(\ref{eq:Gl2-colLie2}),
for the
coloured multiplication extended to the whole dual
$U_{p,q}(\{\alpha,\beta,\gamma,\delta\})$ of~$U^c(gl(1/1))$, are
\begin{eqnarray}
  [\alpha(\lambda), \beta(\mu)]^{\nu} & = & \lambda \omega A^{\mu}_{\nu}\,
         \beta(\nu), \qquad  [\alpha(\lambda), \gamma(\mu)]^{\nu} = \lambda \eta
         A^{\mu}_{\nu}\, \gamma(\nu), \nonumber \\[0cm]
  [\delta(\lambda), \beta(\mu)]^{\nu} & = & \lambda \omega A^{\mu}_{\nu}\,
         \beta(\nu), \qquad  [\delta(\lambda), \gamma(\mu)]^{\nu} = \lambda \eta
         A^{\mu}_{\nu}\, \gamma(\nu), \nonumber \\[0cm]
  [\alpha(\lambda), \delta(\mu)]^{\nu} & = & 0, \qquad  \{\beta(\lambda),
         \gamma(\mu)\}^{\nu} = \beta(\lambda) \nudot \beta(\mu) =
         \gamma(\lambda) \nudot \gamma(\mu) = 0,  \label{eq:Gl1/1-colLie1}
\end{eqnarray}
and
\begin{eqnarray}
  [\alpha(\lambda), \alpha(\mu)]^{\nu} & = & [\delta(\lambda),
           \delta(\mu)]^{\nu} = 0, \nonumber \\
  \alpha(\lambda) \nudot \beta(\mu) & = & A^{\mu}_{\lambda}\, (2\mu)^{-1}
           \left((\lambda+\mu) \alpha(\mu) \nudot \beta(\lambda) + (\lambda-\mu)
           \delta(\mu) \nudot \beta(\lambda)\right), \nonumber \\
  \gamma(\lambda) \nudot \alpha(\mu) & = & A^{\lambda}_{\mu}\, (2\lambda)^{-1}
           \left((\lambda+\mu) \gamma(\mu) \nudot \alpha(\lambda) -
(\lambda-\mu)
           \gamma(\mu) \nudot \delta(\lambda)\right), \nonumber \\
  \alpha(\lambda) \nudot \delta(\mu) & = & \alpha(\mu) \nudot \delta(\lambda)
           + (\lambda-\mu) (\lambda+\mu)^{-1} \left(- \alpha(\mu) \nudot
           \alpha(\lambda) + \delta(\mu) \nudot \delta(\lambda)\right),
\nonumber \\
  \gamma(\lambda) \nudot \beta(\mu) & = & \gamma(\mu) \nudot \beta(\lambda),
           \nonumber \\
  \delta(\lambda) \nudot \beta(\mu) & = & A^{\mu}_{\lambda}\, (2\mu)^{-1}
           \left((\lambda-\mu) \alpha(\mu) \nudot \beta(\lambda) + (\lambda+\mu)
           \delta(\mu) \nudot \beta(\lambda)\right), \nonumber \\
  \gamma(\lambda) \nudot \delta(\mu) & = & A^{\lambda}_{\mu}\, (2\lambda)^{-1}
           \left(- (\lambda-\mu) \gamma(\mu) \nudot \alpha(\lambda) +
           (\lambda+\mu) \gamma(\mu) \nudot \delta(\lambda)\right),
\end{eqnarray}
where
\begin{equation}
  \{x(\lambda), y(\mu)\}^{\nu} \equiv x(\lambda) \nudot y(\mu) + y(\mu) \nudot
  x(\lambda).
\end{equation}
Equation~(\ref{eq:Gl1/1-colLie1}) gives back Eq.~(\ref{eq:Gl1/1-Lie}) whenever
$\lambda$, $\mu$, $\nu \to 1$.\par
%
%
All the results obtained for the dual pair $\left(U^c(gl(1/1)),
Gl^c(1/1)\right)$ are
therefore entirely similar to those derived for $\left(U^c(gl(2)),
Gl^c(2)\right)$ in
Sec.~\ref{sec:gl(2)}.\par
%
%
\section{CONCLUSION}    \label{sec:conclusion}
In the present paper, we extended the notion of dually conjugate Hopf
(super)algebras to the coloured Hopf (super)algebras~${\cal H}^c$,
introduced in~I
and~II. We showed that if the standard Hopf (super)algebras~${\cal H}_q$
that are
the building blocks of~${\cal H}^c$ have Hopf duals~${\cal H}_q^*$, then
the latter
may be used to construct coloured Hopf duals~${\cal H}^{c*}$, endowed with
coloured algebra and antipode maps, but with a standard coalgebraic
structure.\par
%
%
Next, we reviewed the case where the ${\cal H}_q$'s are QUEA's of Lie
(super)algebras~$U_q(g)$, so that the ${\cal H}_q^*$'s are quantum
(super)groups~$G_q$. We extended the Fronsdal and Galindo universal ${\cal
T}$-matrix formalism~\cite{fronsdal93} to the coloured pairs $\left(U^c(g),
G^c\right)$ by introducing coloured universal ${\cal T}$-matrices. We then
proved
that the coloured $RTT$-relations, defining coloured $A(R)$~Hopf
(super)algebras in
Basu-Mallick's approach~\cite{basu}, may be considered as the realization in the
$U_q(g)$~defining representation of a representation-free relation
satisfied by the
coloured universal $\cal R$- and $\cal T$-matrices.\par
%
%
Such results were finally illustrated by constructing two physically-relevant
examples of coloured pairs, corresponding to the two-parameter deformations of
$\bigl(U(gl(2)), Gl(2)\bigr)$, and $\bigl(U(gl(1/1)), Gl(1/1)\bigr)$,
respectively.\par
%
%
In conclusion, we did prove that the formalism developed in I,~II, and the
present
paper, provides an algebraic formulation of the coloured $RTT$-relations, and
establishes a link between the coloured extensions of
Drinfeld-Jimbo~\cite{drinfeld} and Faddeev-Reshetikhin-Takhtajan~\cite{faddeev}
pictures of quantum groups and quantum algebras. Since transfer matrices of
quantum integrable models may be obtained from~$\cal T$ by specialization to
given representations, we do think that the coloured extension of~$\cal T$,
and the
related new algebraic structures introduced here, may have some interesting
applications to such models.\par
%
%
There remain some open questions, which might be interesting topics for future
study. We would like to mention here two of them. The first one consists in
investigating the complementary approach to that considered in the present
paper,
namely trying to transform a set of quantum (super)groups~$G_q$ into a coloured
Hopf (super)algebra~${\cal H}^c$ by defining an appropriate colour group, then
transferring the coloured structure to the duals~$U_q(g)$ to build a
coloured Hopf
dual~${\cal H}^{c*}$. The second problem is to understand the relation, if any,
between the latter and the coloured $U(R)$ Hopf (super)algebras introduced by
Kundu and Basu-Mallick, and defined in terms of coloured
$RLL$-relations~\cite{kundu}.\par
%
%
\newpage
\section*{APPENDIX A: COLOURED MAPS OF $Gl^c(2)$}
\renewcommand{\theequation}{A\arabic{equation}}
\setcounter{equation}{0}
In this Appendix, we list the results obtained for the coloured
multiplication, unit,
and antipode of~$Gl^c(2)$. For their derivation, we used
Eqs.~(\ref{eq:col-maps}),
(\ref{eq:Gl2-alg}), (\ref{eq:Gl2-detcom}), (\ref{eq:Gl2-S}),
(\ref{eq:Gl2-gauss}),
(\ref{eq:Gl2-alphadelta}), (\ref{eq:Gl2-rhoexp}), (\ref{eq:Gl2-rhoainv}), and the
fact that $\bx^{0000}$ is the unit of~$Gl_{p,q}(2)$.\par
%
%
The coloured multiplication is given by
\begin{eqnarray}
  \tm^{\nu}_{p,q,\lambda,\mu}(x \otimes y) & = &
P^{t^{\nu}_{\lambda,\mu}(x,y)/2}\,
         {\cal D}^{(\lambda+\mu-2\nu)/(2\nu)} xy, \label{eq:Gl2-xy} \\
  \tm^{\nu}_{p,q,\lambda,\mu}({\cal D} \otimes x) & = &
         P^{t^{\nu}_{\lambda,\mu}({\cal D},x)/2}\,
         {\cal D}^{(2\lambda+\mu-\nu)/(2\nu)} x, \label{eq:Gl2-Dx} \\
  \tm^{\nu}_{p,q,\lambda,\mu}(x \otimes {\cal D}) & = &
         P^{t^{\nu}_{\lambda,\mu}(x,{\cal D})/2}\,
         {\cal D}^{(\lambda+2\mu-\nu)/(2\nu)} x, \label{eq:Gl2-xD} \\
  \tm^{\nu}_{p,q,\lambda,\mu}({\cal D} \otimes {\cal D}) & = &
         {\cal D}^{(\lambda+\mu)/\nu}, \label{eq:Gl2-DD}
\end{eqnarray}
where
\begin{eqnarray}
  t^{\nu}_{\lambda,\mu}(x,y) & = & 0, \mu-\nu, -\mu+\nu, 0, \lambda+2\mu-3\nu,
         \lambda+3\mu-4\nu, \lambda+\mu-2\nu, \nonumber \\
  & &  \lambda+2\mu-3\nu, -\lambda-2\mu+3\nu, -\lambda-\mu+2\nu,
         -\lambda-3\mu+4\nu, \nonumber \\
  & &  -\lambda-2\mu+3\nu, 0, \mu-\nu, -\mu+\nu, 0,  \label{eq:Gl2-txy} \\
  t^{\nu}_{\lambda,\mu}({\cal D},x) & = & 0, \mu-\nu, -\mu+\nu, 0,
         \label{eq:Gl2-tDx} \\
  t^{\nu}_{\lambda,\mu}(x,{\cal D}) & = & 0, \lambda+4\mu-\nu,
         -\lambda-4\mu+\nu, 0,   \label{eq:Gl2-txD}
\end{eqnarray}
whenever $x$, $y$ run over $\{a,b,c,d\}$, and are listed in lexicographical
order.\par
%
%
For the coloured unit and antipode, the results read
\begin{equation}
  \tiota^{\nu}_{p,q}(1_k) = \tiota_{p^{(\nu)},q^{(\nu)}}(1_k),
\end{equation}
and
\begin{eqnarray}
  \tS^{\nu}_{p,q,\mu}(a) & = & {\cal D}^{-(\mu+\nu)/(2\nu)} d, \qquad
        \tS^{\nu}_{p,q,\mu}(b) = - Q^{-1} \left(P^{\nu}
        {\cal D}\right)^{-(\mu+\nu)/(2\nu)} b, \\
  \tS^{\nu}_{p,q,\mu}(c) & = & - Q \left(P^{-\nu}
        {\cal D}\right)^{-(\mu+\nu)/(2\nu)} c, \qquad \tS^{\nu}_{p,q,\mu}(d) =
        {\cal D}^{-(\mu+\nu)/(2\nu)} a, \\
  \tS^{\nu}_{p,q,\mu}({\cal D}) & = & {\cal D}^{-\mu/\nu},
\end{eqnarray}
respectively.\par
%
%
\newpage
\section*{APPENDIX B: COLOURED MAPS OF $Gl^c(1/1)$}
\label{sec:append-Gl1/1}
\renewcommand{\theequation}{B\arabic{equation}}
\setcounter{equation}{0}
In this Appendix, we list the results obtained for the coloured
multiplication, unit,
and antipode of~$Gl^c(1/1)$. For their derivation, we used
Eqs.~(\ref{eq:col-maps}),
(\ref{eq:Gl2-gauss}), (\ref{eq:Gl1/1-alg}), (\ref{eq:Gl1/1-S}),
(\ref{eq:Gl1/1-alphadelta}), (\ref{eq:Gl1/1-rhoexp}),
(\ref{eq:Gl1/1-rhoainv}), and
the fact that $\bx^{0000}$ is the unit of~$Gl_{p,q}(1/1)$.\par
%
%
The coloured multiplication is given by
\begin{eqnarray}
  \tm^{\nu}_{p,q,\lambda,\mu}(x \otimes y) & = & C^{\nu}_{\lambda,\mu}(x,y)
         {\cal D}^{(2\nu-\lambda-\mu)/(2\nu)} a^{(\lambda+\mu-2\nu)/\nu} xy,
         \label{eq:Gl1/1-xy} \\
  \tm^{\nu}_{p,q,\lambda,\mu}(a \otimes d) & = &
         {\cal D}^{(2\nu-\lambda-\mu)/(2\nu)} a^{(\lambda+\mu-2\nu)/\nu} \left\{
         ad - p^{\mu} q^{-\nu} \left(A^{\nu-\mu}_{\nu}\right)^2 cb\right\}, \\
  \tm^{\nu}_{p,q,\lambda,\mu}(d \otimes a) & = &
         {\cal D}^{(2\nu-\lambda-\mu)/(2\nu)} a^{(\lambda+\mu-2\nu)/\nu}
\biggl\{
         da + p^{\mu} \biggl(q^{-\nu} \left(A^{\nu-\mu}_{\nu}\right)^2
\nonumber \\
         & & \mbox{} - p^{\lambda}
         \left(A^{2\nu-\lambda-\mu}_{\nu}\right)^2\biggr) cb\biggr\}, \\
  \tm^{\nu}_{p,q,\lambda,\mu}(d \otimes d) & = &
         {\cal D}^{(2\nu-\lambda-\mu)/(2\nu)} a^{(\lambda+\mu-2\nu)/\nu} \left\{
         d^2 - p^{\lambda+\mu} \left(A^{2\nu-\lambda-\mu}_{\nu}\right)^2
         ca^{-1}db\right\},  \label{eq:Gl1/1-dd} \\
  \tm^{\nu}_{p,q,\lambda,\mu}({\cal D} \otimes x) & = &
         C^{\nu}_{\lambda,\mu}({\cal D},x) {\cal D}^{(3\nu-\mu)/(2\nu)}
         a^{(\mu-\nu)/\nu} x, \label{eq:Gl1/1-Dx} \\
  \tm^{\nu}_{p,q,\lambda,\mu}({\cal D} \otimes d) & = &
         {\cal D}^{(3\nu-\mu)/(2\nu)} a^{(\mu-\nu)/\nu} \left\{d - p^{\mu}
         \left(A^{\nu-\mu}_{\nu}\right)^2 ca^{-1}b\right\}, \\
  \tm^{\nu}_{p,q,\lambda,\mu}(x \otimes {\cal D}) & = &
         C^{\nu}_{\lambda,\mu}(x,{\cal D}) {\cal D}^{(3\nu-\lambda)/(2\nu)}
         a^{(\lambda-\nu)/\nu} x, \label{eq:Gl1/1-xD} \\
  \tm^{\nu}_{p,q,\lambda,\mu}(d \otimes {\cal D}) & = &
         {\cal D}^{(3\nu-\lambda)/(2\nu)} a^{(\lambda-\nu)/\nu} \left\{d -
p^{\lambda}
         \left(A^{\nu-\lambda}_{\nu}\right)^2 ca^{-1}b\right\}, \\
  \tm^{\nu}_{p,q,\lambda,\mu}({\cal D} \otimes {\cal D}) & = & {\cal D}^2.
         \label{eq:Gl1/1-DD}
\end{eqnarray}
In Eq.~(\ref{eq:Gl1/1-xy}), $x$, $y$ run over $\{a,b,c,d\}$ with the
exceptions of
$(x,y) = (a,d)$, $(d,a)$, $(d,d)$, and
\begin{eqnarray}
  C^{\nu}_{\lambda,\mu}(x,y) & = & 1, A^{\mu}_{\nu}, q^{\mu-\nu} A^{\mu}_{\nu},
         p^{\mu-\nu} A^{\lambda}_{\nu}, 1, (pq)^{\mu-\nu} A^{\lambda}_{\nu}
         A^{\mu}_{\nu}, p^{\mu-\nu} A^{\lambda}_{\nu}, q^{\lambda+\mu-2\nu}
         A^{\lambda}_{\nu},  \nonumber \\
  & & q^{\lambda+\mu-2\nu} A^{\lambda}_{\nu} A^{\mu}_{\nu}, 1,
         q^{\lambda+\mu-2\nu} A^{\lambda}_{\nu}, A^{\mu}_{\nu}, q^{\mu-\nu}
         A^{\mu}_{\nu}, \label{eq:Gl1/1-Cxy}
\end{eqnarray}
for $(x,y)$ listed in lexicographical order. Similarly, in
Eqs.~(\ref{eq:Gl1/1-Dx})
and~(\ref{eq:Gl1/1-xD}), $x$  runs over $\{a,b,c\}$, and
\begin{equation}
  C^{\nu}_{\lambda,\mu}({\cal D},x) = 1, A^{\mu}_{\nu}, q^{\mu-\nu}
A^{\mu}_{\nu},
  \qquad C^{\nu}_{\lambda,\mu}(x,{\cal D}) = 1, A^{\lambda}_{\nu},
q^{\lambda-\nu}
  A^{\lambda}_{\nu}.  \label{eq:Gl1/1-CDx}
\end{equation}
\par
%
%
For the coloured unit and antipode, the results read
\begin{equation}
  \tiota^{\nu}_{p,q}(1_k) = \tiota_{p^{\nu},q^{\nu}}(1_k),
\end{equation}
and
\begin{eqnarray}
  \tS^{\nu}_{p,q,\mu}(a) & = & {\cal D}^{(\mu-\nu)/(2\nu)} \left\{a^{-\mu/\nu}
        - p^{-\mu} q^{\mu} \left(A^{\mu}_{\nu}\right)^2 c
a^{-(\mu+\nu)/\nu} \hd^{-1}
        b\right\}, \\
  \tS^{\nu}_{p,q,\mu}(b) & = & - p^{-\mu} A^{\mu}_{\nu}\, a^{-(\mu+\nu)/(2\nu)}
        \hd^{-(\mu+\nu)/(2\nu)}\, b, \\
  \tS^{\nu}_{p,q,\mu}(c) & = & - q^{\mu} A^{\mu}_{\nu}\, c\,
a^{-(\mu+\nu)/(2\nu)}
        \hd^{-(\mu+\nu)/(2\nu)}, \\
  \tS^{\nu}_{p,q,\mu}(d) & = & {\cal D}^{-(\mu-\nu)/(2\nu)} \hd^{-\mu/\nu}, \\
  \tS^{\nu}_{p,q,\mu}({\cal D}) & = & {\cal D}^{-1},
\end{eqnarray}
respectively.\par
%
%


\newpage
\begin{thebibliography}{99}
%
\bibitem{cq96a} C. Quesne, ``Coloured quantum universal enveloping
algebras,'' J.
Math. Phys. (in press).
%
\bibitem{cq96b} C. Quesne, ``Coloured Hopf algebras,'' in {\it Proc. XXI
Int. Coll. on
Group Theoretical Methods in Physics, Goslar, Germany, July 15--20, 1996} (in
press), q-alg/9705019.
%
\bibitem{bazhanov} V. V. Bazhanov and Yu. G. Stroganov, Theor. Math. Phys.
{\bf 62},
253 (1985); L. Hlavat\'y, J. Phys. A {\bf 20}, 1661 (1987); R. J. Baxter,
J. H. Perk,
and H. Au-Yang, Phys. Lett. A {\bf 128}, 138 (1988); B. S. Shastry, J.
Stat. Phys. {\bf
50}, 57 (1988).
%
\bibitem{footnote1} The coloured YBE, considered here and in~I, should not be
confused with the colour YBE~\cite{mcanally} that arises when extending the
graded
YBE~\cite{chaichian} to more general gradings than that determined by
$\Z_2$~\cite{rittenberg}. As a consequence, our coloured Hopf algebras are
distinct
from the Hopf colour algebras~\cite{mcanally}, generalizing Hopf
superalgebras~\cite{chaichian} to such more general gradings.
%
\bibitem{drinfeld} V. G. Drinfeld, in {\it Proc. Int. Congress of
Mathematicians (Berkeley, CA, 1986)\/}, edited by A. M. Gleason (AMS,
Providence,
RI, 1987), p. 798; M. Jimbo, Lett. Math. Phys. {\bf 10}, 63 (1985); {\bf
11}, 247
(1986).
%
\bibitem{cq97} C. Quesne, ``Coloured Hopf algebras,'' in {\it Proc. Workshop on
Special Functions and Differential Equations, Chennai, India, January
13--24, 1997}
(in press), q-alg/9705022.
%
\bibitem{chaichian} M. Chaichian and P. Kulish, Phys. Lett. B {\bf 234}, 72
(1990);
W. B. Schmidke, S. P. Volos, and B. Zumino, Z. Phys. C {\bf 48}, 249 (1990).
%
\bibitem{faddeev} L. Faddeev, N. Reshetikhin, and L. Takhtajan, in {\it
Algebraic
Analysis\/}, Vol.~1, edited by M. Kashiwara and T. Kawai (Academic, New York,
1988) p.~129; in {\it Braid Group, Knot Theory and Statistical
Mechanics\/}, edited
by C. N. Yang and M. L. Ge (World Scientific, Singapore, 1989) p.~97.
%
\bibitem{ohtsuki} T. Ohtsuki, J. Knot Theor. Its Rami. {\bf 2}, 211 (1993).
%
\bibitem{bonatsos} D. Bonatsos, P. Kolokotronis, C. Daskaloyannis, A. Ludu,
and C.
Quesne, Czech. J. Phys. {\bf 46}, 1189 (1996); D. Bonatsos, C. Daskaloyannis, P.
Kolokotronis, A. Ludu, and C. Quesne, J. Math. Phys. {\bf 38}, 369 (1997).
%
\bibitem{kundu} A. Kundu and B. Basu-Mallick, J. Phys. A {\bf 27}, 3091
(1994); B.
Basu-Mallick, Mod. Phys. Lett. A {\bf 9}, 2733 (1994).
%
\bibitem{basu} B. Basu-Mallick, Int. J. Mod. Phys. A {\bf 10}, 2851 (1995).
%
\bibitem{majid} S. Majid, Int. J. Mod. Phys. A {\bf 5}, 1 (1990); V. Chari
and A.
Pressley, {\it A Guide to Quantum Groups} (Cambridge U.P., Cambridge, 1994).
%
\bibitem{fronsdal93} C. Fronsdal and A. Galindo, Lett. Math. Phys. {\bf 27}, 59
(1993).
%
\bibitem{fronsdal94} C. Fronsdal and A. Galindo, Contemp. Math. {\bf 175}, 73
(1994); C. Fronsdal, in {\it Noncompact Lie Groups and Some of Their
Applications
(San Antonio, TX, 1993)}, NATO Adv. Sci. Inst. Ser. C Math. Phys. Sci., 429
(Kluwer,
Dordrecht, 1994) p.~423.
%
\bibitem{bonechi} F. Bonechi, E. Celeghini, R. Giachetti, C. M. Pere\~na,
E. Sorace, and M.
Tarlini, J. Phys. A {\bf 27}, 1307 (1994).
%
\bibitem{morozov} A. Morozov and L. Vinet, ``Free-field representation of group
element for simple quantum groups,'' Universit\'e de Montr\'eal preprint
CRM-2202,
hep-th/9409093 (1994).
%
\bibitem{chakra96a} R. Chakrabarti and R. Jagannathan, Z. Phys. C {\bf 72}, 519
(1996).
%
\bibitem{finkelstein} R. J. Finkelstein, Lett. Math. Phys. {\bf 29}, 75 (1993).
%
\bibitem{jaga} R. Jagannathan and J. Van der Jeugt, J. Phys. A {\bf 28}, 2819
(1995); J. Van der Jeugt and R. Jagannathan, Czech. J. Phys. {\bf 46}, 269
(1996).
%
\bibitem{chakra96b} R. Chakrabarti and R. Jagannathan, Lett. Math. Phys.
{\bf 37},
191 (1996).
%
\bibitem{schirrmacher} A. Schirrmacher, J. Wess, and B. Zumino, Z. Phys. C
{\bf 49},
317 (1991).
%
\bibitem{dobrev} V. K. Dobrev, J. Math. Phys. {\bf 33}, 3419 (1992).
%
\bibitem{burdik92a} \v C. Burd\'\i k and P. Hellinger, J. Phys. A {\bf 25}, L629
(1992).
%
\bibitem{chakra94} R. Chakrabarti and R. Jagannathan, J. Phys. A {\bf 27}, 2023
(1994).
%
\bibitem{footnote2} In~I, as well as in Secs.~\ref{sec:duals}
and~\ref{sec:QUEA} of
the present paper, we used the generic symbol~$q^{\nu}$, where $\nu$
behaved as a
contravariant index, to denote the parameters of the final algebra under colour
group transformations. Since in this Section, and in
Sec.~\ref{sec:gl(1/1)}, powers
of the parameters make their appearance, to avoid confusion we slightly
modify our
notation~$q^{\nu}$ into~$q^{(\nu)}$.
%
\bibitem{footnote3} It should be noted that the images of $a$, $b$, $c$,~$d$
under~$\rho^{\nu}$ or~$\rho_{\nu}$ belong to the larger enveloping algebra
$U_{p,q}(\{\alpha,\beta,\gamma,\delta\})$.
%
\bibitem{dabrowski} L. Dabrowski and L. Wang, Phys. Lett. B {\bf 266}, 51
(1991).
%
\bibitem{chakra91} R. Chakrabarti and R. Jagannathan, J. Phys. A {\bf 24}, 5683
(1991).
%
\bibitem{bednar} M. Bedn\'a\v r, \v C. Burd\'\i k, M. Couture, and L.
Hlavat\'y, J. Phys.
A {\bf 25}, L341 (1992); \v C. Burd\'\i k and R. Tom\'a\v sek, Lett. Math.
Phys. {\bf
26}, 97 (1992).
%
\bibitem{burdik92b} \v C. Burd\'\i k and P. Hellinger, J. Phys. A {\bf 25},
L1023
(1992).
%
\bibitem{chakra95} R. Chakrabarti and R. Jagannathan, Z. Phys. C {\bf 66}, 523
(1995).
%
\bibitem{kauffman} L. H. Kauffman and H. Saleur, Commun. Math. Phys. {\bf 141},
293 (1991); L Rozanskey and H. Saleur, Nucl. Phys. B {\bf 376}, 461 (1992).
%
\bibitem{hinrichsen} H. Hinrichsen and V. Rittenberg, Phys. Lett. B {\bf
275}, 350
(1992).
%
\bibitem{footnote4} Although we use the same symbol to denote the duals
of~$U_{p,q}(gl(2))$ and~$U_{p,q}(gl(1/1))$, there is of course no relation
between
them.
%
\bibitem{mcanally} D. S. McAnally, in {\it Proc. Yamada Conf. XL, XX Int.
Coll. on
Group Theoretical Methods in Physics, Toyonaka, Japan, July 4--9, 1994\/},
edited
by A. Arima, T. Eguchi, and N. Nakanishi (World Scientific, Singapore,
1995) p.~339.

%
\bibitem{rittenberg} V. Rittenberg and D. Wyler, Nucl. Phys. B {\bf 139}, 189
(1978); J. Math. Phys. {\bf 19}, 2193 (1978); J. Lukierski and V.
Rittenberg, Phys.
Rev. D {\bf 18}, 385 (1978); M. Scheunert, J. Math. Phys. {\bf 20}, 712 (1979).

\end{thebibliography}
\end{document}